\documentclass[12pt]{article}
\usepackage{amsmath}
\usepackage{amssymb}
\usepackage{amscd}
\usepackage{psfrag}
\usepackage{graphicx}
\newcommand{\pr}{{\rm pr}}
\newcommand{\T}{{\rm T}}
\newcommand{\reg}{{\rm reg}}
\newcommand{\FF}{{\bf F}}
\newcommand{\ind}{{\bf 1}}
\newcommand{\PR}{{\rm PR}}
\newcommand{\Diag}{{\rm Diag}}
\newcommand{\W}{{\cal W}}
\newcommand{\Disc}{{\rm D}}
\newcommand{\Ball}{{\rm B}}
\newcommand{\V}{{\cal V}}
\newcommand{\PH}{{\rm PH}}
\newcommand{\PSH}{{\rm PSH}}
\renewcommand{\P}{\mathbb{P}}
\newcommand{\C}{\mathbb{C}}
\newcommand{\B}{{\cal B}}
\newcommand{\R}{\mathbb{R}}

\newcommand{\id}{{\rm id}}

\renewcommand{\H}{{\cal H}}
\newcommand{\G}{{\cal G}}

\newcommand{\PC}{{\rm PC}}
\newcommand{\aire}{{\rm aire}}

\newcommand{\ddc}{{\rm dd^c}}
\renewcommand{\d}{{\rm d}}

\newcommand{\cad}{{\it c.-\`a-d. }}
\newcommand{\ie}{{\it i.e. }}

\newcommand{\Lone}{{{\rm L}^1}}
\newcommand{\Ltwo}{{{\rm L}^2}}
\newcommand{\Linfty}{{{\rm L}^\infty}}

\newcommand{\Ltwoloc}{{{\rm L}^2_{\rm loc}}}

\newcommand{\diam}{{\rm diam}}

\newcommand{\F}{{\cal F}}
\newcommand{\A}{{\cal A}}
\newcommand{\supp}{{\rm supp}}
\newcommand{\E}{{\cal E}}

\renewcommand{\O}{{\rm O}}
\renewcommand{\o}{{\rm o}}

\newtheorem{theoreme}{Th\'eor\`eme}[section]
\newtheorem{proposition}[theoreme]{Proposition}
\newtheorem{corollaire}[theoreme]{Corollaire}
\newtheorem{lemme}[theoreme]{Lemme}

\newtheorem{exemples}[theoreme]{Exemples}
\newtheorem{remarque}[theoreme]{Remarque}

\newenvironment{preuve}{\begin{trivlist} \item[]{\bf D\'emonstration.}}
{\par\hfill $\square$\end{trivlist}}
\title{Distribution des pr\'eimages et 
des points p\'eriodiques d'une correspondance polynomiale} 
\author{Tien-Cuong Dinh}
\date{}
\begin{document}
\maketitle
\hfill \`a Madame L\^e H\^ong S\^am
\\
\begin{abstract} We construct an equilibrium measure $\mu$ for a
polynomial correspondence $F$ of Lojasiewicz exponent
$l>1$. We then show
that $\mu$ can be built as the distribution of preimages of
a generic point and that the repelling periodic points are
equidistributed on the support of $\mu$. Using this results, we
will give a characterization of infinite
uniqueness sets for polynomials.  
\end{abstract}
\\
{\it Mots cl\'es:} Correspondance, Mesure d'\'equilibre, 
Ensemble exceptionnel, Point\break
p\'eriodique, Ensemble d'unicit\'e
\\
\normalsize
\section{Introduction}
Un compact $K$ de $\C$ est un {\it ensemble d'unicit\'e} 
si pour tout couple de polyn\^omes non constants $f$ et $g$ la relation
$f^{-1}(K)=g^{-1}(K)$ implique $f=g$.
Ostrovskii, Pakovitch et Zaidenberg \cite{Ostrovskii} ont
montr\'e que si $f$ et $g$ sont deux polyn\^omes de
m\^eme degr\'e v\'erifiant $f^{-1}(K)=g^{-1}(K)$ pour un compact
$K$ de cardinal au moins deux,
alors il existe une rotation $R$
pr\'eservant $K$ telle que $f=R\circ g$. Dans \cite{Dinh} nous avons
d\'etermin\'e les polyn\^omes $f$, $g$ et les compacts $K$ de
capacit\'e logarithmique positive v\'erifiant
$f^{-1}(K)=g^{-1}(K)$. Des probl\`emes analogues 
pour les fonctions enti\`eres
ou m\'eromorphes  ont \'et\'e \'etudi\'es par 
Nevanlinna \cite{Nevanlinna}, Gross-Yang \cite{GrossYang},
Shiffman \cite{Shiffman}... Ces auteurs utilisent des m\'ethodes vari\'ees.
\par
Ici, \`a partir de la relation $f^{-1}(K)=g^{-1}(K)$ on d\'eduit que $g\circ
f^{-1}(K)=K$ et par cons\'equent
$$(g\circ f^{-1})\circ \cdots
\circ (g\circ f^{-1})(K)=K,$$
ce qui permet de se ramener \`a l'\'etude dynamique de la fonction
multivalu\'ee $F:=f\circ g^{-1}$ qui est en fait une correspondance
polynomiale. Les propri\'et\'es dynamiques que nous allons
\'etudier
permettent de caract\'eriser les ensembles d'unicit\'e
infinis (voir le corollaire 5.2).
\par 
Dans \cite{ClozelUllmo1, ClozelUllmo2}, Clozel et Ullmo
\'etudient les correspondances holomorphes sur les surfaces de
Riemann et sur les domaines sym\'etriques. Ils en donnent des applications
en arithm\'etique. Ils montrent que les correspondances
modulaires sur une courbe holomorphe hyperbolique sont celles qui 
pr\'eservent la forme volume $\Omega$
associ\'ee
\`a la m\'etrique de Kobayashi. Ils en d\'eduisent que les
correspondances, qui commutent avec une correspondance modulaire ext\'erieure,
sont modulaires car ces correspondances, elles aussi, doivent
pr\'eserver $\Omega$.
\par
Depuis les travaux de Julia \cite{Julia}, Fatou \cite{Fatou},
Ritt \cite{Ritt}, Eremenko \cite{Eremenko} (voir aussi
\cite{DinhSibony3}), on sait que
si deux endomorphismes de $\P^k$ commutent, les objets dynamiques,
qui leur sont associ\'es,
sont fortement li\'es. L'\'etude de ces objets permet de
d\'eterminer ou de caract\'eriser ces endomorphismes.
Dans le cadre des correspondances holomorphes,
une \'etude dynamique
devrait permettre de comprendre les commutateurs (voir le
corollaire 2.9). 
\par
Nous renvoyons le lecteur \`a \cite{BedfordSmillie, Fornaess,
FornaessSibony,RussakovskiShiffman,Sibony},
pour les aspects fondamentaux de la th\'eorie
d'it\'eration des applications holomorphes et m\'eromorphes de
$\P^k$.  
Pour les
endomorphismes holomorphes de $\P^k$ ou pour
les automorphismes de H\'enon de $\C^2$ par exemple, on sait
construire des mesures invariantes,
m\'elangeantes qui maximisent
l'entropie. Ces mesures d'\'equilibre sont obtenues
comme intersections de courants invariants positifs
ferm\'es de bidegr\'e $(1,1)$.
\par
Briend-Duval \cite{BriendDuval1,BriendDuval2}
ont montr\'e que la mesure d'\'equilibre de tout endomorphisme
holomorphe de degr\'e $d\geq 2$ de $\P^k$
est limite de masses de Dirac port\'ees par
les points p\'eriodiques
r\'epulsifs. C'est aussi la limite de masses de Dirac port\'ees par les
pr\'eimages de tout point $z$ n'appartenant
pas \`a un ensemble exceptionnel alg\'ebrique $\E$.
Ant\'erieurement, Forn\ae ss-Sibony \cite{FornaessSibony}
avaient montr\'e que $\E$ est pluripolaire. En dimension 1, ces
r\'esultats ont \'et\'e d\'emontr\'es par Brolin pour les
polyn\^omes \cite{Brolin},
Lyubich \cite{Lyubich} et Freire-Lop\`es-Ma\~ n\'e
\cite{FreireLopesMane} pour les fractions rationnelles.
Notons ici que l'\'etude des
endomorphismes holomorphes de $\P^k$, peut se ramener \`a
l'\'etude des endomorphismes polynomiaux. Il suffit de consid\'erer le relev\'e de ces
applications \`a $\C^{k+1}$. Dans \cite{DinhSibony1}, l'auteur
et Sibony ont construit pour les applications d'allure polynomiale
une mesure invariante d'entropie maximale
et g\'en\'eralis\'e les th\'eor\`emes de Briend-Duval pour de
grandes familles de telles applications (en particulier pour les
applications polynomiales dont l'exposant de Lojasiewicz est
sup\'erieur \`a 1).
\par
Dans le pr\'esent
travail, nous allons g\'en\'eraliser les r\'esultats de
Briend-Duval aux
correspondances polynomiales. Notre article s'organise de
la mani\`ere suivante. Au paragraphe 2 nous d\'efinissons les
correspondances polynomiales sur $\C^k$
et leurs exposants de Lojasiewicz \`a
l'infini. Nous construisons la mesure d'\'equilibre $\mu$
associ\'ee \`a une correspondance polynomiale $F$ d'exposant de
Lojasiewicz $l>1$. Cette mesure est $F^*$-invariante,
``m\'elangeante'' \`a vitesse
exponentielle et ne charge pas les ensembles pluripolaires.
Nous montrons aussi que toute correspondance polynomiale, qui
commute avec $F$, pr\'eserve la mesure d'\'equilibre $\mu$ de $F$
(voir le corollaire 2.9).
La construction de $\mu$ suit une m\'ethode donn\'ee dans
\cite{DinhSibony1} ({\it m\'ethode par r\'esolution de $\ddc$}); 
elle est aussi valable pour les
correspondances {\it d'allure polynomiale} ou pour les
it\'erations al\'eatoires (voir aussi \cite{DinhSibony4}).
Pour certaines correspondances, on peut construire un courant
invariant $T$
positif ferm\'e de bidegr\'e $(1,1)$. Mais il est peu
probable que la mesure $T^k$ (m\^eme lorsqu'elle est bien
d\'efinie) soit invariante quand $k\geq 2$. 
\par
Dans le troisi\`eme
paragraphe,
en adaptant les m\'ethodes
de Lyubich \cite{Lyubich} et Briend-Duval \cite{BriendDuval1,
BriendDuval2} (voir aussi \cite{DinhSibony1, Dinh2}),
nous construisons, pour les petites boules centr\'ees en
un point g\'en\'erique, beaucoup de branches inverses dont on contr\^ole 
la taille.
La mesure $\mu$ r\'efl\`ete la distribution des
pr\'eimages de tout point $z$ qui n'appartient pas \`a un ensemble
exceptionnel $\E$. En collaboration avec Charles Favre, nous
montrons que $\E$ est l'orbite positive de $\E_0$ o\`u $\E_0$ est
le plus grand sous-ensemble alg\'ebrique propre
de $\C^k$ invariant par $F^{-1}$. Le cas des applications
polynomiales d'exposant de Lojasiewicz $l>1$ est trait\'e dans
\cite{DinhSibony1, DinhSibony2} (voir aussi \cite{Guedj,DinhSibony4}). 
On obtient alors que l'ensemble
$\E$ est {\it alg\'ebrique}.
\par
Dans le quatri\`eme paragraphe, nous montrons en particulier  
que les points p\'eriodiques r\'eguliers r\'epulsifs de
$F$ sont denses et \'equidistribu\'es sur le support de $\mu$.
\par
Les r\'esultats obtenus sont encore valables dans un cadre plus g\'en\'eral.
Afin de simplifier les notations,
nous pr\'ef\'erons
nous limiter au cas de l'espace complexe $\C^k$.
Dans \cite{DinhSibony4, Dinh2}, nous \'etendons cette \'etude aux it\'erations 
al\'eatoires des correspondances sur les vari\'et\'es k\"ahl\'eriennes compactes. 
\par
Une interpr\'etation g\'eom\'etrique des r\'esultats obtenus est
donn\'ee \`a la fin du paragraphe 4 (voir les corollaires 4.7,
4.8). Cette vision g\'eom\'etrique
nous semble int\'eressante m\^eme pour les endomorphismes
holomorphes de $\P^k$ et les automorphismes de H\'enon de $\C^2$.
Dans le dernier paragraphe, nous appliquons les r\'esultats
obtenus pour d\'eterminer les ensembles d'unicit\'e infinis, pour
les polyn\^omes d'une variables.
\par
Signalons un
travail r\'ecent de Claire Voisin \cite{Voisin} dans lequel
elle \'etudie la non-hyperbolicit\'e de vari\'et\'es
projectives en utilisant des correspondances (voir aussi l'exemple
3.12). Du point de vue
dynamique, les correspondances consid\'er\'ees par Clozel-Ullmo
et Claire Voisin sont plus proches des automorphismes holomorphes
tandis que celles \'etudi\'ees dans le pr\'esent article sont
plut\^ot proches des endomorphismes holomorphes de $\P^k$.
Un outil que nous avons d\'evelopp\'e  r\'ecemment avec Sibony
\cite{DinhSibony5} permet
d'\'etudier les correspondances de type automorphisme.
\par
Dans la suite, $\Ball(z,r)$ et $\overline \Ball(z,r)$
d\'esignent la boule ouverte et la boule ferm\'ee
de centre $z$ et de rayon $r$. Les disques,
les boules 
et le diam\`etre $\diam(.)$ d'un ensemble sont d\'efinis ou
mesur\'es en m\'etrique
euclidienne. L'aire $\aire(.)$ d'un disque, la masse
$\|.\|$ d'un courant, les normes $\Ltwo$ et ${\cal C}^2$ d'une
fonction sont
mesur\'es en m\'etrique de Fubini-Study. La notation $\delta_z$
d\'esigne la masse de Dirac en $z$, $\ind_S$ d\'esigne la
fonction indicatrice de $S$. Les pr\'eimages d'un point $z$ de
$F$ sont aussi les images de $z$ par la correspondance $\overline F$
adjointe \`a $F$. Nous pr\'ef\'erons parler de pr\'eimages plut\^ot que
d'images afin
que les applications polynomiales soient couvertes par notre
\'etude.
\par
\ 
\\
{\bf Remerciements.} Je remercie Charles Favre et Nessim
Sibony dont les nombreuses
remarques ont permis d'am\'eliorer la r\'edaction de cet article.  
\section{Correspondances polynomiales}
Soit $X$ une vari\'et\'e complexe de dimension $k\geq 1$. Notons
$\pi_1$, $\pi_2$ les projections canoniques de $X\times X$ dans
$X$. On appelle {\it $k$-cha\^{\i}ne holomorphe}
de $X\times X$ toute
combinaison finie $Y:=\sum n_i Y_i$ o\`u les $Y_i$ sont des
sous-ensembles analytiques irr\'eductibles de dimension $k$, 
deux \`a deux distincts, de $X\times X$ et
o\`u les $n_i$ sont des entiers relatifs non nuls.
On dira que $Y$ est {\it positive}
si les $n_i$ sont
positifs.
D'apr\`es un th\'eor\`eme de Lelong,
une $k$-cha\^{\i}ne holomorphe positive $Y$ d\'efinit par
int\'egration un courant positif ferm\'e $[Y]$ de bidimension
$(k,k)$ de $X\times X$. 
Notons $|Y|:=\cup Y_i$ le support de $Y$ et 
$\overline Y:=\sum n_i \overline Y_i$ o\`u $\overline Y_i$ 
est le sym\'etrique de $Y_i$ par rapport \`a la diagonale 
de $X\times X$, \ie l'image de $Y_i$ par l'application
$(x,y)\mapsto (y,x)$.  
\par 
Une {\it correspondance holomorphe} de degr\'e topologique
$(d_1,d_2)$
sur $X$ est la donn\'ee d'une $k$-cha\^{\i}ne holomorphe positive 
$Y$ de dimension $k$ de $X\times X$ telle que la restriction de
$\pi_i$ \`a $Y$ d\'efinisse une application propre de degr\'e $d_i$
pour $i=1,2$. Plus
pr\'ecis\'ement, pour tout $z\in X$ la fibre $Y\cap
\pi_i^{-1}(z)$ contient exactement $d_i$ points compt\'es avec
multiplicit\'es.
Il est clair que si $(x,y)\in |Y|$ et si $x$ tend vers l'infini alors
$y$ tend aussi vers l'infini et r\'eciproquement.
On peut identifier cette correspondance \`a 
la fonction multivalu\'ee
$F:=(\pi_{2|Y})\circ(\pi_{1|Y})^{-1}$. Le
terme {\it correspondance} d\'esignera $F$.
On dira que $Y$ est le {\it graphe} de $F$.
On utilisera
souvent 
la d\'ecomposition $Y=\sum Y_i^*$ dans laquelle
chaque $Y_i$ est r\'ep\'et\'e
$n_i$ fois afin d'\'eviter de parler de multiplicit\'es.
La correspondance $\overline F$ associ\'ee
\`a $\overline Y$
est appel\'ee {\it correspondance adjointe} de $F$.
\par
Soit $F'$ une autre correspondance de degr\'e topologique
$(d_1',d_2')$
associ\'ee \`a une $k$-cha\^{\i}ne holomorphe positive
$Z=\sum Z_j^*$. La composition
$F'\circ F$ est celle associ\'ee
au produit fibr\'e $Y\times_X Z:=\sum (Y^*_i\times_X Z^*_j)$
o\`u 
\begin{eqnarray*}
Y^*_i\times_X Z^*_j & := & \big\{(x,z)\in X\times X
\mbox{ tel qu'il existe } \\
& & \ \ y\in X \mbox{ v\'erifiant }
(x,y)\in Y_i^* \mbox{ et } (y,z)\in Z_j^* \big\}.
\end{eqnarray*}
Le produit $Y^*_i\times_X Z^*_j$ est, en g\'en\'eral, une
$k$-cha\^{\i}ne holomorphe qui n'est pas toujours irr\'eductible.
La composition $F'\circ F$
est une correspondance de degr\'e topologique $(d_1d_1',d_2d_2')$.
On notera $F^n$ la correspondance
$F\circ\cdots \circ F$ ($n$ fois).
\par
Dans le pr\'esent travail, nous consid\'erons le cas o\`u
$X$ est l'espace euclidien 
$\C^k$ et les composantes $Y_i$ de $Y$ 
sont des sous-ensembles alg\'ebriques
de $\C^k\times\C^k$. 
On dira qu'une telle correspondance $F$ est 
{\it polynomiale (propre)}. Il existe une constante $l>0$ 
telle que pour tout $(x,y)\in |Y|$ suffisamment grand
on ait $|y|\geq c|x|^l$
o\`u $c>0$ est une constante.
Si $F$ est un endomorphisme polynomial, Ploski \cite{Ploski} a
montr\'e qu'il existe une constante maximale $l>0$ qui v\'erifie
la propri\'et\'e ci-dessus. Sa preuve est aussi valable
pour les correspondances polynomiales. 
Cette constante $l$ 
est appel\'ee {\it exposant de Lojasiewicz} de $F$.
Dans la suite, on suppose que 
$l>1$. On v\'erifie que dans ce cas $d_1$ est
strictement plus petit que $d_2$ (on peut prouver ceci en
utilisant l'argument donn\'e dans la proposition 4.1).
Notons $z$ les coordonn\'ees euclidiennes de $\C^k$ et 
$\omega:=\frac{1}{2}\ddc\log(1+\|z\|^2)$
la forme de Fubini-Study de $\P^k$.
Soit $A_0>1$ une constante assez grande que nous allons choisir
dans le lemme 2.5. 
Fixons un nombre $R_0>0$ assez grand tel que 
$|y|>A_0|x|$ pour tout $(x,y)\in |Y|$ v\'erifiant $|y|\geq R_0$.
\par
On pose $F^{-1}:=(\pi_{1|Y})\circ (\pi_{2|Y})^{-1}$,
$F^*:=(\pi_{1|Y})_*(\pi_{2|Y})^*$ et
$F_*:=(\overline F)^*=(\pi_{2|Y})_*(\pi_{1|Y})^*$.
Les "applications" $F$ et $F^{-1}$ agissent sur les
sous-ensembles de $\C^k$,  
les points
de la fibre $F^{-1}(z)$ de $F$
sont compt\'es avec multiplicit\'es. L'op\'erateur $F_*$ 
agit sur les fonctions continues ou
plurisousharmoniques (p.s.h.) et
sur les courants positifs ferm\'es de bidegr\'e $(1,1)$ de
$\C^k$. L'op\'erateur $F^*$ agit
sur les mesures positives. Plus pr\'ecis\'ement, si $\varphi$ est
une fonction continue ou p.s.h. sur $\C^k$, on pose
$$F_*\varphi:=\sum_{w\in F^{-1}(z)}\varphi(w).$$
C'est une fonction continue ou p.s.h. sur $\C^k$.
Rappelons que $\varphi$ est {\it p.s.h.}
si elle est localement int\'egrable, 
semi-continue sup\'erieurement (s.c.s)
et si $\ddc\varphi\geq 0$ au sens des courants.
Elle est {\it pluriharmonique}
si elle est continue et $\ddc\varphi=0$. Observons que si $\varphi$
est pluriharmonique, $F_*\varphi$ l'est aussi.
Si $T$ est un
courant positif ferm\'e de bidegr\'e $(1,1)$ sur $\C^k$, il existe une
fonction p.s.h. $\varphi$, unique \`a une fonction
pluriharmonique pr\`es, telle que $\ddc\varphi=T$. On d\'efinit
alors $F_*T:=\ddc F_*\varphi$ (voir \cite{Meo}).
Pour une mesure
positive $\nu$ \`a support compact
sur $\C^k$, on d\'efinit $F^*\nu$ par
$$\langle F^*\nu,\varphi\rangle:=\langle\nu,F_*\varphi \rangle
\ \mbox{ pour } \varphi \mbox{ continue sur }\C^k.$$
\par
Nous dirons qu'une fonction p.s.h. $\varphi$ sur $\C^k$ est {\it \`a 
croissance logarithmique} s'il existe une constante $A>0$ telle que 
$\varphi-A\log(1+\|z\|^2)$ soit born\'ee sup\'erieurement.
On dit qu'une mesure positive est {\it PB} si 
elle int\`egre les fonctions p.s.h. 
\`a croissance logarithmique \cite{DinhSibony1}.
Dans le cas de dimension $k=1$, si $\mu$ est une mesure positive sur $\C$, on 
peut \'ecrire $\mu=\ddc u-\alpha$ avec $u$ une fonction $\Lone$ et $\alpha$ une
forme lisse sur $\P^1$. 
On a montr\'e \cite{DinhSibony1} que $\mu$ est PB 
si et seulement si son Potentiel $u$ est Born\'e. Ceci justifie la terminologie
choisie. 
D'apr\`es le th\'eor\`eme de Josefson \cite[Theorem
5.2.4]{Klimek}, pour tout ensemble pluripolaire $E$, il existe
une fonction p.s.h. \`a croissance logarithmique $\varphi$ telle
que $\varphi=-\infty$ sur $E$. Par cons\'equent, 
les mesures PB ne chargent pas les ensembles pluripolaires.
En particulier, le support d'une mesure PB est parfait, \ie ne
contient pas de point isol\'e.
\begin{theoreme} Soit $F$ une correspondance polynomiale 
de degr\'e topologique
$(d_1,d_2)$ sur $\C^k$, d'exposant
de Lojasiewicz $l>1$. 
Soient $\nu_n$ des mesures de 
probabilit\'e de support uniform\'ement born\'e dans 
$\C^k$. Supposons que 
$\nu_n=h_n\omega^k$ o\`u $h_n$ est une fonction v\'erifiant 
$\|h_n\|_{\Ltwo}=\o(l^n)$.
Alors la suite de mesures
$d_2^{-n}(F^n)^*\nu_n$ converge vers une mesure de
probabilit\'e $\mu$, \`a support compact, 
ind\'ependante de la suite 
$(\nu_n)$. De plus, cette mesure est PB
et v\'erifie la relation de
$F^*$-invariance: $F^*\mu=d_2\mu$.
\end{theoreme}
On dit que $\mu$ est la {\it mesure d'\'equilibre} de $F$. En
g\'en\'eral, elle n'est pas invariante par $F$, \ie $F_*\mu\not=d_1\mu$.
Pour montrer le th\'eor\`eme 2.1,
l'id\'ee est de tester une fonction p.s.h. $\varphi$. On a 
$$d_2^{-n}\langle (F^n)^*(\nu_n),\varphi\rangle =  \langle \nu_n, d_2^{-n} 
(F^n)_*\varphi \rangle.$$
Nous allons montrer que $d_2^{-n} (F^n)_*\varphi$ tend dans $\Ltwoloc$ \`a vitesse 
$\O(l^{-n})$ vers
une constant $c_\varphi$. Le th\'eor\`eme en d\'ecoule.
\par
Fixons une boule $V:=\Ball(0,R)$ de rayon $R>R_0$ qui contient les supports 
des $\nu_n$. 
Posons $U:=F^{-1}(V)$. Alors $\overline U$ est contenu dans la boule 
$U':=\Ball(0,R/A_0)$. Fixons aussi $V':=\Ball(0,R')$
avec $R'>R$ tel que 
$F^{-1}(V')\subset U'$. 
Soit $\varphi$ une fonction sur $V$ ou sur $\C^k$.
Soit $\Lambda:=d_2^{-1}F_*$
l'op\'erateur de Perron-Frobenius associ\'e \`a $F$. On a par
d\'efinition 
$$\Lambda\varphi(z) =d_2^{-1}\sum_{w\in F^{-1}(z)}\varphi(w)$$
o\`u les points de $F^{-1}(z)$ sont compt\'es avec multiplicit\'es.
Posons $\varphi_n:=\Lambda^n \varphi$. 
Cette fonction est p.s.h. ou continue si $\varphi$ 
l'est. 
Nous aurons besoin des lemmes suivants.
\begin{lemme} Soit $T$ un courant positif ferm\'e de
bidegr\'e $(1,1)$
et de masse $1$ sur $\P^k$. Alors la masse de
$\Lambda^n T$ dans $\C^k$ est 
plus petite ou \'egale \`a $l^{-n}$.
\end{lemme}
\begin{preuve} Dans $\C^k$, 
on peut \'ecrire $T=\ddc\varphi$ o\`u $\varphi$ est une fonction
p.s.h. telle que $\varphi-\frac{1}{2}\log(1+\|z\|^2)$
soit born\'ee sup\'erieurement. On a 
$$l^n\Lambda^nT=\ddc (l^n\Lambda^n\varphi)=\ddc 
(l^n\varphi_n).$$
Comme l'exposant de Lojasiewicz de $F$ est \'egal \`a $l>1$, quitte
\`a effectuer un changement lin\'eaire de coordonn\'ees, on peut
supposer que $\|z'\|^l\leq \|z\|$ lorsque $\|z\|$ est assez grand
et $z'\in F^{-1}(z)$.  
La fonction $l^n\varphi_n-\frac{1}{2}\log(1+\|z\|^2)$ est donc
born\'ee sup\'erieurement car les valeurs de $\varphi_n(z)$ sont
obtenues comme la moyenne des valeurs de $\varphi$ sur $F^{-n}(z)$. 
D'apr\`es un lemme de comparaison (voir par exemple
\cite{Sibony}, \cite[proposition 5.4]{DinhSibony2}), ceci implique que
\begin{eqnarray*}
\|l^nd_2^{-n}(F^n)_*T\|_{\C^k} & = & \int_{\C^k}\ddc
(l^n\varphi_n)\wedge \omega^{k-1}\\
& \leq & \frac{1}{2}\int_{\C^k}\ddc \log(1+\|z\|^2)\wedge \omega^{k-1}
=\int_{\C^k} \omega^k=1.
\end{eqnarray*} 
\end{preuve}
\begin{lemme} Soit $\varphi$ une fonction
p.s.h. sur $V$ telle que 
$\varphi_n:=\Lambda^n\varphi$ ne tende pas uniform\'ement vers
$-\infty$ (en particulier si 
$\varphi$ est localement born\'ee). Alors
$\varphi_n$ tend vers une constante
$c_\varphi$ dans l'espace $\Ltwoloc(\C^k)$. 
\end{lemme}
\begin{preuve} 
Observons que $\varphi_n$ est d\'efinie dans la
boule $\Ball(0,A_0^nR)$. On 
montre que $\varphi_n$ tend vers une constante
$c_\varphi$ dans $\Ltwoloc(V)$.
Pour la convergence dans $\Ltwoloc(\C^k)$ il
suffit de remplacer $V$ par 
des boules suffisamment grandes.
\par 
Soit $\psi$ la r\'egularis\'ee s.c.s. de  
la fonction $(\limsup\varphi_n)$ sur $V$. 
C'est une fonction p.s.h. car $\varphi_n$ ne tend pas 
uniform\'ement vers $-\infty$.  Comme 
$F^{-1}(V)= U$, on a $\sup_V \varphi_{n+1}\leq
\sup_U\varphi_n$. En effet, les valeurs
de $\varphi_{n+1}$ dans $V$ sont obtenues comme la moyenne
de valeurs de 
$\varphi_n$ dans $U$.
Ceci implique que $\sup_V\psi\leq \sup_U\psi$ car 
$\psi$ est aussi \'egale \`a la
r\'egularis\'ee s.c.s. de la fonction $(\limsup\varphi_{n+1})$.
Le principe du maximum entra\^{\i}ne que $\psi$ est une 
constante $c_\varphi$.
\par
Soit $(\varphi_{n_i})$ une sous-suite convergente
dans $\Ltwoloc(V)$ vers 
une fonction p.s.h. $h$. Montrons que $h=c_\varphi$.
On a $h\leq c_\varphi$.
Si $h\not=c_\varphi$, le principe du maximum implique que
$h\leq c_\varphi-2\epsilon$ sur
$U$ o\`u $\epsilon>0$
est une constante. D'apr\`es le lemme de Hartogs \cite[Theorem
2.6.4]{Hormander},
on a $\varphi_{n_i}
\leq c_\varphi-\epsilon$ pour $i$ assez grand. Par cons\'equent,
$\varphi_n\leq c_\varphi-\epsilon$ sur $V$ pour tout $n>n_i$.
Ceci contredit
le fait que la r\'egularis\'ee s.c.s. de
$(\limsup\varphi_n)$ est \'egale \`a $c_\varphi$. 
On a montr\'e que $\varphi_n$ tend 
vers $c_\varphi$ dans $\Ltwoloc(V)$.
\par 
\end{preuve}
\par
Le lemme pr\'ec\'edent permet de construire
{\it la mesure d'\'equilibre} 
$\mu$. 
Si $\varphi$ est une fonction de classe ${\cal C}^2$ \`a
support compact, la suite 
$\varphi_n$ converge aussi vers une constante $c_\varphi$ dans 
$\Ltwoloc(\C^k)$
car cette fonction $\varphi$ s'\'ecrit comme diff\'erence de deux 
fonctions p.s.h.
Si $\Omega$ est une forme volume lisse \`a support dans $V$ telle que 
$\int\Omega=1$, on a d'apr\`es le lemme 2.3
$$\lim_{n\rightarrow \infty}
\langle d_2^{-n}(F^n)^*\Omega,\varphi\rangle=
\lim_{n\rightarrow\infty}
\langle\Omega,\Lambda^n\varphi\rangle= 
c_\varphi.$$
Par cons\'equent, $d_2^{-n}(F^n)^*\Omega$ tend faiblement vers une  
mesure de probabilit\'e $\mu$ d\'efinie par
$\langle\mu,\varphi\rangle:=c_\varphi$ pour 
$\varphi$ de classe ${\cal C}^2$. Il est clair
que $\mu$ est port\'ee par
$\overline U$ et v\'erifie la relation de $F^*$-invariance.
Elle ne d\'epend pas de $\Omega$.
Observons ici que la fonction
p.s.h. $\varphi$ v\'erifie 
l'hypoth\`ese du lemme 2.3 si et seulement
si elle est $\mu$-int\'egrable.
\begin{lemme} Soit $\varphi$ une fonction p.s.h. sur $V$. 
Si $\varphi$ n'est pas $\mu$-int\'egrable, 
la suite de fonctions $\varphi_n:=\Lambda^n\varphi$ converge 
uniform\'ement vers $-\infty$. Si $\varphi$ est $\mu$-int\'egrable, $(\varphi_n)$ 
converge dans $\Ltwoloc(\C^k)$ vers la constante $c_\varphi:=\int\varphi\d\mu$. 
\end{lemme}
\begin{preuve} Comme $\mu$ est $F^*$-invariante, on a $\langle \mu,\varphi\rangle = 
\langle\mu,\varphi_n\rangle$ pour tout $n\geq 0$. 
D'apr\`es le lemme 2.3, il suffit de traiter le cas o\`u 
la suite $(\varphi_n)$ converge dans $\Ltwoloc(\C^k)$ 
vers une constante $c_\varphi$. Montrons que 
$c_\varphi=\int\varphi\d\mu$. 
\par
D'apr\`es le lemme de Hartogs (voir le lemme 2.3), 
on a $\limsup\varphi_n=c_\varphi$. On d\'eduit de la 
relation $\langle \mu,\varphi\rangle=\langle \mu,\varphi_n\rangle$
que $\langle \mu,\varphi\rangle\leq c_\varphi$. 
D'autre part, la semi-continuit\'e
sup\'erieure de $\varphi$ et la d\'efinition de $\mu$ impliquent que
$$\langle\mu,\varphi\rangle\geq \limsup \langle d_t^{-n} 
(F^n)^*\Omega,\varphi \rangle = \limsup\langle\Omega,\varphi_n\rangle =c_\varphi.$$
La preuve du lemme est achev\'ee.
\end{preuve}
\begin{lemme} Soit $\varphi$ une fonction
pluriharmonique sur $V$. Alors on a pour $A_0$ assez grand
$$\|\varphi_n-c_\varphi\|_{\Linfty(V)}\leq 2^{-n}l^{-n}
\|\varphi-c_\varphi\|_{\Linfty(V)}.$$ 
\end{lemme}
\begin{preuve} Quitte \`a remplacer $\varphi$
par $\varphi-c_\varphi$, 
on peut supposer que $c_\varphi=\langle\mu,\varphi\rangle=0$.
Le support de $\mu$ \'etant contenu dans la boule
$U'=\Ball(0,R/A_0)$, l'\'egalit\'e $\langle\mu,\varphi\rangle=0$
implique que $\varphi$ 
doit s'annuler en un point de
cette boule. D'apr\`es un lemme du type Schwarz, on a
$$\|\varphi\|_{\Linfty(U)}\leq (2l)^{-1}\|\varphi\|_{\Linfty(V)}$$
lorsque $A_0$ est assez grand.
Par cons\'equent, $\|\varphi_1\|_{\Linfty(V)}\leq (2l)^{-1}\|\varphi
\|_{\Linfty(V)}$. Ceci implique le lemme. Notons que pour montrer le lemme du type 
Schwarz, on peut consid\'erer la famille normale des fonctions harmoniques
$\psi$ sur le disque unit\'e $\Delta\subset\C$ s'annulant en un point
de $\Delta(0,1/A_0)$ et v\'erifiant $\psi\leq 1$.
\end{preuve}
\par
Notons $\PSH(V)$ le c\^one des fonctions p.s.h. sur
$V$. Rappelons que les boules $U'$ et $V'$
sont fix\'ees au d\'ebut de la d\'emonstration du th\'eor\`eme
2.1. 
\begin{lemme}
L'op\'erateur $\Lambda:\PSH(U')\cap\Ltwo(U')\longrightarrow
\PSH(V)\cap\Ltwo(V)$ est born\'e dans le sens o\`u
il existe une constante $c>0$ telle que
$\|\Lambda\varphi\|_{\Ltwo(V)}\leq c\|\varphi\|_{\Ltwo(U')}$ pour
toute fonction $\varphi\in\PSH(U')\cap \Ltwo(U')$.
\end{lemme}
\begin{preuve} Soit $\varphi^{(n)}$ une suite de fonctions p.s.h.
sur $U'$ telle que $\|\varphi^{(n)}\|_{\Ltwo(U')}\leq 1$. Il faut
montrer que la suite $\Lambda\varphi^{(n)}$ est
born\'ee dans $\Ltwo(V)$.
\par
Quitte \`a extraire une sous-suite, on peut supposer que la suite
$\varphi^{(n)}$ converge dans $\Ltwoloc(U')$
vers une fonction p.s.h. $\varphi$. On en d\'eduit que
la suite de fonctions p.s.h.
$\Lambda\varphi^{(n)}$ converge vers $\Lambda\varphi$
dans $\Ltwoloc(V')$.
En particulier, elle converge vers $\Lambda\varphi$ dans $\Ltwo(V)$.  
\end{preuve}
\par
Notons $\PH(V)$ l'espace des fonctions pluriharmoniques sur $V$.
Posons $H:=\PH(V)\cap \Ltwo(V)$ et $H^\perp$
son orthogonal dans $\Ltwo(V)$.
Posons aussi $H^{\perp *}:=\PSH(V)\cap H^\perp$.
Le sous-espace $H$ est invariant par $\Lambda$ car si $\varphi$ est 
pluriharmonique sur $V$, $\Lambda\varphi$ l'est sur un voisinage de $\overline V$.
Pour toute fonction
$\varphi\in \PSH(V')$, on a la d\'ecomposition unique
$\varphi=u+v$ avec $u\in H$ et $v\in H^{\perp *}$.
La fonction $v$ est le potentiel dans $V$ du courant $\ddc\varphi$
dont la norme
$\Ltwo$ est minimale. 
\par
D'apr\`es le lemme 2.6, l'op\'erateur 
$\Lambda: \PSH(V)\cap\Ltwo(V)
\longrightarrow \PSH(V)\cap\Ltwo(V)$ est born\'e.
Par cons\'equent, il existe des applications
lin\'eaires born\'ees
$\Lambda_1:H\longrightarrow H$,
$\Lambda_2:H^{\perp *}\longrightarrow H$ et
$\Lambda_3:H^{\perp *}\longrightarrow H^{\perp *}$
telles que $\Lambda\varphi=\Lambda_1 u+\Lambda_2v+\Lambda_3v$.
On a $\Lambda_1=\Lambda_{|H}$ et
$\Lambda_2=\pr_H\circ\Lambda_{|H^{\perp *}}$ et
$\Lambda_3=\pr_{H^\perp}\circ\Lambda_{|H^{\perp *}}$ o\`u
$\pr_{|H}$ et $\pr_{|H^\perp}$
d\'esignent les projections orthogonales de $\Ltwo(V)$
sur $H$ et sur $H^\perp$.
On a $\ddc\Lambda_3^n\varphi=\ddc\varphi_n$.
\begin{proposition} La mesure $\mu$ est PB. 
De plus, il existe une constante $c>0$
telle que pour toute fonction p.s.h. $\varphi$ avec 
$\|\varphi\|_{\Ltwo(V)}\leq A$ et 
$\varphi-\frac{1}{2}A\log(1+\|z\|^2)$ born\'ee sup\'erieurement, on ait
$\|\Lambda^n\varphi-c_\varphi\|_{\Ltwo(V)}
\leq cA l^{-n}$ pour tout $n\geq 1$.
\end{proposition}
\begin{preuve} Par homoth\'etie, il suffit de
consid\'erer une fonction 
$\varphi$ p.s.h. avec
$\|\varphi\|_{\Ltwo(V)}\leq 1$ et
$\varphi-\frac{1}{2}\log(1+\|z\|^2)$ born\'ee sup\'erieurement.
La famille de telles fonctions est compacte. 
Les constantes $c_i$ et $c$ que nous allons utiliser 
sont positives et ind\'ependantes de $\varphi$. 
D'apr\`es le lemme 2.2, on a $\|\ddc\varphi_n\|\leq l^{-n}$. Par cons\'equent,
il existe une fonction p.s.h. $\psi_n$ telle que $\psi_n-\frac{1}{2}l^{-n}
\log(1+\|z\|^2)$ soit born\'ee sup\'erieurement, $\|\psi_n\|_{\Ltwo}\leq c_1l^{-n}$
et $\ddc\psi_n=\ddc\varphi_n$. 
Il s'agit ici la r\'esolution de $\ddc$ sur $\P^k$.
On en d\'eduit que
$\|\Lambda_3^n v\|_{\Ltwo(V)}\leq\|\psi_n\|_\Ltwo\leq c_1l^{-n}$.
Posons
$$b:=\int u\d\mu,\ \ b_n:=\int\Lambda_2\Lambda_3^n v\d\mu
\ \mbox{ et }\ 
s_n:=b+b_1+\cdots+b_{n-1}.$$
La fonction pluriharmonique $u$ v\'erifie
$\|u\|_{\Ltwo(V)}\leq \|\varphi\|_{\Ltwo(V)}\leq 1$, on en d\'eduit par la formule
de la moyenne, que
$\|u\|_{\Linfty(U)}\leq c_2$.
Comme $\Lambda_2$ est born\'e, on a
$$\|\Lambda_2\Lambda_3^n v\|_{\Ltwo(V)}\leq c_3\|\Lambda_3^n
v\|_{\Ltwo(V)} \leq c_1c_3 l^{-n}.$$
La formule de la moyenne
appliqu\'ee \`a la fonction pluriharmonique
$\Lambda_2\Lambda_3^n v$ 
implique que $|b_n|\leq c_4 l^{-n}$
et donc la suite $(s_n)$ est convergente.
Utilisant les in\'egalit\'es obtenues et
le lemme 2.5 et l'estimation de $b_j$, on obtient en
d\'eveloppant $\Lambda^n$,
\begin{eqnarray*}
\|\Lambda^n\varphi-s_n\|_{\Ltwo(V)}
& = & \|\Lambda_1^nu+\Lambda_1^{n-1}\Lambda_2 v
+\Lambda_1^{n-2}\Lambda_2\Lambda_3 v + \cdots +
\Lambda_1\Lambda_2\Lambda_3^{n-2} v +\\
& & + \Lambda_2\Lambda_3^{n-1} v
-s_n\|_{\Ltwo(V)} +\|\Lambda_3^n v\|_{\Ltwo(V)}\\
& \leq &
\|\Lambda_1^nu-b\|_{\Ltwo(V)}
+\|\Lambda_1^{n-1}\Lambda_2 v-b_1\|_{\Ltwo(V)}+\\
& & +\|\Lambda_1^{n-2}\Lambda_2\Lambda_3 v -b_2\| + \cdots 
+\|\Lambda_1\Lambda_2\Lambda_3^{n-2} v-b_{n-2}\| +\\
& & + \|\Lambda_2\Lambda_3^{n-1}v -b_{n-1}\|_{\Ltwo(V)}+
\|\Lambda_3^nv\|_{\Ltwo(V)}\\
&\leq  & c_5\big[(2l)^{-n}+(2l)^{-n+1}+ (2l)^{-n+2}l^{-1}+\\
& & +\cdots + (2l)^{-1}l^{-n+2}+l^{-n+1}+l^{-n}\big]\\
&\leq & c_6 l^{-n}.
\end{eqnarray*}
On en d\'eduit que $\Lambda^n\varphi$ converge vers la constante
$\lim s_n=b+\sum b_i$. Par cons\'equent, $\varphi$ est
$\mu$-int\'egrable et $c_\varphi=b+\sum b_i$. On a aussi
$|c_\varphi-s_n|\leq c_7l^{-n}$. Ceci implique que
$\|\Lambda^n\varphi-c_\varphi\|_{\Ltwo(V)}\leq cl^{-n}$.
\end{preuve}
\begin{corollaire} 
Il existe une constante
  $c>0$ telle que $\sup_V(\varphi_n-c_\varphi)\leq cAl^{-n}$ 
pour tout $n\geq 1$ et 
toute fonction $\varphi$ p.s.h. sur $\C^k$ avec 
$\|\varphi\|_{\Ltwo(V')}\leq A$ et
$\varphi-\frac{1}{2}A\log(1+\|z\|^2)$ 
born\'ee sup\'erieurement. 
\end{corollaire}
\begin{preuve} Il suffit d'appliquer 
la proposition 2.7 en rempla\c cant $V$ par une boule $V_1$ telle que 
$\overline V\subset V_1$ et $\overline V_1\subset V'$, puis 
d'utiliser l'in\'egalit\'e de la sous-moyenne pour les fonctions
p.s.h.
\end{preuve}
{\bf Fin de la d\'emonstration du th\'eor\`eme 2.1.} Soit
$\varphi$ une fonction p.s.h. \`a croissance logarithmique.
D'apr\`es l'in\'egalit\'e de Cauchy-Schwarz, on a
\begin{eqnarray*}
|\langle d_2^{-n}(F^n)^*\nu_n-\mu,\varphi\rangle| & = &
|\langle\nu_n,\Lambda^n\varphi\rangle -c_\varphi|\\
& = & |\langle h_n\omega^k,\Lambda^n\varphi-c_\varphi\rangle| \\
& \leq & \|h_n\|_{\Ltwo(V)}
\|\varphi_n -c_\varphi\|_{\Ltwo(V)}.
\end{eqnarray*}
D'apr\`es la proposition 2.7, la derni\`ere expression tend vers
0 car $\|h_n\|_{\Ltwo(V)}=\o(l^n)$.  
On en d\'eduit que $d_2^{-n}(F^n)^*\nu_n$ tend faiblement
vers $\mu$.
\par
\hfill $\square$
\par
\begin{corollaire} Soit $F$ une correspondance comme au
th\'eor\`eme 2.1. Si $G$ est une autre correspondance polynomiale
de degr\'e topologique $(d_1',d_2')$, d'exposant de Lojasiewicz
quelconque et
v\'erifiant $F\circ G= G\circ F$ alors $G^*\mu=d_2'\mu$.
\end{corollaire}
\begin{preuve} Soit $\nu$ une forme volume lisse \`a support compact
sur $\C^k$ telle
que $\int_{\C^k}\nu=1$. D'apr\`es le th\'eor\`eme 2.1,
on a $\lim d_2^{-n}(F^n)^*\nu=\mu$ et
$\lim d_2^{-n} (F^n)^*(G^*\nu)=d_2'\mu$. Du fait que $F$ et $G$
commutent, on a
$$\lim_{n\rightarrow\infty} d_2^{-n}(F^n)^*(G^*\nu)=
\lim_{n\rightarrow\infty} d_2^{-n} G^*(F^n)^*\nu=G^*\mu.$$
On en d\'eduit que $G^*\mu=d_2'\mu$.
\end{preuve}
\begin{theoreme} Soit $F$ une correspondance
comme au th\'eor\`eme 2.1. 
Alors la
mesure $\mu$ est $F^*$-m\'elangeante
d'ordre exponentiel. Plus
pr\'ecis\'ement, il existe une constante $c>0$ telle que
pour toute fonction $\varphi$ de classe ${\cal
C}^2$ et toute fonction $\psi$ born\'ee dans $\C^k$, on ait
$|I_n(\varphi,\psi)| \leq cl^{-n}
\|\varphi\|_{{\cal C}^2}\|\psi\|_{\Linfty}$
o\`u 
$$I_n(\varphi,\psi):=\int (\Lambda^n\varphi)\psi\d\mu -
\left(\int\varphi\d\mu\right) \left(\int\psi\d\mu\right).$$
\end{theoreme}
\begin{preuve} 
Consid\'erons d'abord le cas o\`u $\varphi$ est p.s.h. avec
$\|\varphi\|_{\Ltwo(V')}\leq A$ et 
$\varphi-\frac{1}{2}A\log(1+\|z\|^2)$ born\'ee sup\'erieurement. 
Montrons que $I_n(\varphi,\psi)\leq 
cAl^{-n}\|\psi\|_{\Linfty}$.
Observons que si $\alpha$ est 
une constante, on a $I_n(\varphi,\psi+\alpha)=
I_n(\varphi,\psi)$ car $\mu$ est 
$F^*$-invariante. Il suffit donc d'examiner
le cas o\`u $\psi$ est 
positive. On a d'apr\`es le corollaire 2.8
$$I_n(\varphi,\psi)=\int(\Lambda^n\varphi -c_\varphi)\psi \d\mu
\leq \int cAl^{-n} \psi \d\mu \leq cAl^{-n}\|\psi\|_{\Linfty}.$$
De la m\^eme mani\`ere, on montre que $-I_n(\varphi,\psi)=
I_n(\varphi,-\psi)\leq  cAl^{-n}\|\psi\|_{\Linfty}$.
Par cons\'equent, on a
$|I_n(\varphi,\psi)|\leq  cAl^{-n}\|\psi\|_{\Linfty}$.
\par
Pour le cas g\'en\'eral, on peut supposer que
$\varphi$ est \`a support 
compact dans $V$.
Elle s'\'ecrit comme diff\'erence de deux fonctions p.s.h. 
v\'erifiant les conditions ci-dessus. La constante
$A$ est de l'ordre de
$\|\varphi\|_{{\cal C}^2}$. On est ramen\'e au cas pr\'ec\'edent. 
\end{preuve}
\section{Equidistribution des pr\'eimages}
Dans ce paragraphe, nous \'etudions la distribution
des pr\'eimages de
$F$. Pour tout
point $z\in\C^k$ et tout $n\geq 0$,
posons $\mu^z_n:=d_2^{-n} (F^n)^*\delta_z$. Notons $\E$
l'ensemble des points $z$ tels que la suite de mesures
$(\mu^z_n)$ ne tend pas vers $\mu$. C'est {\it l'ensemble
exceptionnel} de $F$.  
Nous avons la proposition suivante.
\begin{proposition} Soit $F$ une correspondance comme
au th\'eor\`eme 2.1. 
Alors $\E$ est pluripolaire.
\end{proposition}
\par
Observons que puisque la mesure $\mu$ est PB, elle ne charge pas 
les ensembles pluripolaires. La proposition entra\^{\i}ne que 
$\mu(\E)=0$ et $\mu^z_n$ tend
faiblement vers $\mu$ pour $\mu$-presque tout point $z\in\C^k$.
On montrera plus loin que $\E$ est une r\'eunion finie ou
d\'enombrable d'ensembles alg\'ebriques.
\begin{preuve} 
Consid\'erons la fonction
  p.s.h. $\varphi:=\log(1+\|z\|^2)$.
C'est une fonction strictement p.s.h.
Posons $c_\varphi:=\int \varphi\d\mu$ et 
$$\Phi:=\sum_{n=0}^\infty (\Lambda^n\varphi-c_\varphi)=
\sum_{n=0}^\infty (\varphi_n -c_\varphi).$$ 
On peut appliquer la proposition 2.7 et le corollaire
2.8 \`a une boule $V$ 
arbitrairement grande. On d\'eduit que $\Phi$ est une fonction
p.s.h. sur $\C^k$. Notons $\E^*$ l'ensemble o\`u
$\Phi$ prend la valeur
$-\infty$. C'est un ensemble pluripolaire. Il suffit de montrer
que $\E\subset \E^*$. 
\par
Soit $z\not\in \E^*$.
Du fait que $\Phi(z)$ est finie, la suite $\varphi_n(z)$ tend vers
$c_\varphi$.
Soit $\psi$ une fonction r\'eelle 
${\cal C}^2$ \`a support compact. Pour montrer
que $\mu^z_n$ tend vers 
$\mu$, il suffit
de montrer que $\psi_n(z)$ tend vers $c_\psi:=\int\psi\d\mu$. Fixons 
$\epsilon>0$ suffisamment petit tel que
$\varphi^{\pm}:=\varphi\pm\epsilon\psi$ soit p.s.h. 
Un tel $\epsilon$ existe car $\varphi$
est strictement p.s.h.
D'apr\`es le lemme 2.7, les suites de
fonctions $\varphi^\pm_n:=\Lambda^n\varphi^\pm$
convergent vers les constantes
$c^\pm:=\int\varphi^\pm \d \mu= c_\varphi\pm\epsilon c_\psi$ dans
$\Ltwoloc(\C^k)$. 
Comme on l'a d\'ej\`a vu dans la preuve du lemme 2.3, 
$\limsup\varphi^\pm_n(z)$ est au plus \'egal \`a $c^\pm$. 
Or $\lim\varphi_n(z)=c_\varphi$ car $z\not\in \E^*$.
On en d\'eduit que
$\limsup\pm\psi_n(z)\leq \pm c_\psi$, donc $\lim\psi_n(z)=c_\psi$. 
\end{preuve}
\par
Dans la suite, nous allons d\'ecrire plus pr\'ecis\'ement la
distribution des pr\'eimages de $F$. Nous allons
en fait construire des branches inverses holomorphes d\'efinies sur 
des disques et des boules.   
\par
Rappelons que $Y=\sum Y_i^*$ est le graphe de $F$.
Soit $K_0$ un sous-ensemble connexe de $\C^k$.
On appelle {\it branche
inverse \underline{r\'eguli\`ere} d'ordre $n$}
de $K_0$ (voir Figure 1) toute suite $\B$ 
$$K_{-n},(\widehat{K}_{-n},i_n), K_{-n+1},
(\widehat{K}_{-n+1},i_{n-1}),
\ldots,K_{-1}, (\widehat{K}_{-1},i_1), K_0$$
v\'erifiant les propri\'et\'es suivantes
\begin{enumerate}
\item[(i)] Les ensembles $K_{-m}\subset\C^k$ et
$\widehat{K}_{-m}\subset Y^*_{i_m}$ sont connexes;
\item[(ii)] $\pi_1$ d\'efinit une bijection de
$\widehat{K}_{-m}$ dans $K_{-m}$ et
$\pi_2$ d\'efinit une bijection de
$\widehat{K}_{-m}$ dans $K_{-m+1}$ 
pour tout $1\leq m\leq n$.
\end{enumerate}
\begin{figure}
\begin{center}
\psfrag{A}{$K_0$}
\psfrag{B}{$K_{-1}$}
\psfrag{C}{$K_{-2}$}
\psfrag{E}{$\widehat K_{-1}$}
\psfrag{F}{$\widehat K_{-2}$}
\psfrag{K}{$Y^*_1$}
\psfrag{Y}{$Y^*_2$}
\psfrag{K}{$Y^*_1$}
\psfrag{H}{$\C^k$}
\includegraphics{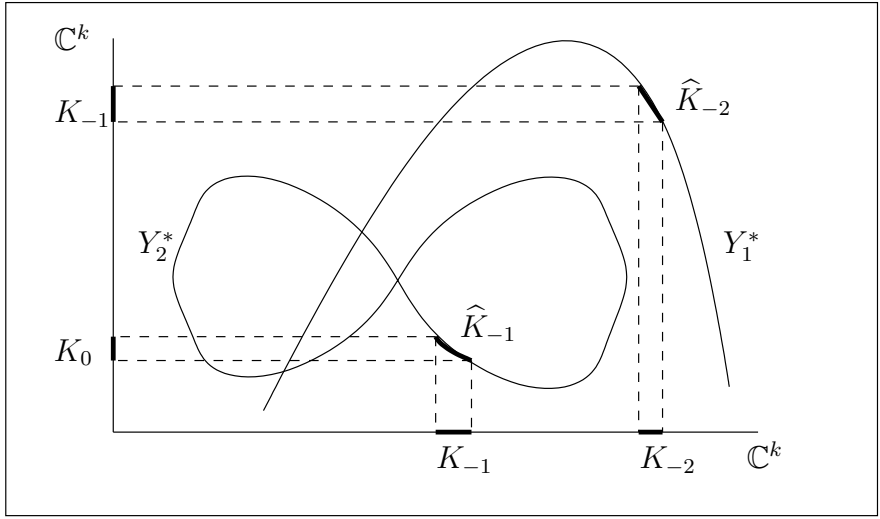}
\caption{Branche inverse r\'eguli\`ere.}
\end{center}
\end{figure}
Si $K_0$ n'est pas un ouvert (par exemple si $K_0$ est un point),
on exige que les applications
ci-dessus d\'efinissent des bijections holomorphes
entre un voisinage de $\widehat K_{-m}$ et ses images.
Puisque les $Y_i^*$ ne sont pas n\'ecessairement distincts, 
les indices $i_m$ permettent de compter les branches inverses
r\'eguli\`eres avec multiplicit\'e.
Il y a au plus
$d_2^n$ branches inverses r\'eguli\`eres d'ordre $n$ de $K_0$. 
En pratique, $K_0$ sera un point, un disque holomorphe,
une famille de disques centr\'es en un point ou
une boule holomorphe.
Notons $F^{-m}_\B$ l'application
$$(\pi_{1|\widehat K_{-m}})\circ (\pi_{2|\widehat K_{-m}})^{-1}
\circ \cdots \circ
(\pi_{1|\widehat K_{-1}})\circ (\pi_{2|\widehat K_{-1}})^{-1}$$
pour $1\leq m\leq n$.
C'est une application holomorphe bijective de $K_0$ dans $K_{-m}$.
Soit $(a_m)$ une suite de nombres r\'eels tendant vers 0.
On dira que
$\B$ est {\it de taille $(a_m)$} si
le diam\`etre de l'ensemble $K_{-m}$ est au plus \'egal \`a $a_m$
pour $0\leq m\leq n$. Soient $z_0\in K_0$ et 
$\B_0$ une branche inverse r\'eguli\`ere d'ordre $n$ du point $z_0$ 
donn\'ee par la suite 
$$z_{-n},(\widehat{z}_{-n},i_n), z_{-n+1},
(\widehat{z}_{-n+1},i_{n-1}),
\ldots,z_{-1}, (\widehat{z}_{-1},i_1), z_0$$
(les points $\widehat{z}_{-m}$ appartiennent \`a $Y^*_{i_m}$ et sont 
de coordonn\'ees $(z_{-m+1},z_{-m})$).
On dira que la branche $\B$ est {\it accroch\'ee \`a  
la branche $\B_0$} si on a 
$\widehat{z}_{-m}\in \widehat{K}_{-m}$ pour tout $1\leq m\leq n$
(les indices $i_n$ dans $\B_0$ et $\B$ sont identiques).
\par
Nous soulignons ici que deux
branches inverses r\'eguli\`eres
$\B$ et $\B'$ de $K_0$ sont consid\'er\'ees comme \'egales si et
seulement si on
a $K_{-m}=K_{-m}'$, $\widehat K_{-m}=\widehat K_{-m}'$ ainsi que
$i_m=i_m'$ pour tout $m$ o\`u $K'_{-m}$, $\widehat K'_{-m}$
et $i_m'$ sont des ensembles et des indices associ\'es \`a $\B'$.
Nous allons donner plus loin une notion de branche inverse plus
souple qui permet de prouver que $\E$ est une r\'eunion
d'ensembles alg\'ebriques. 
Pour les branches \underline{r\'eguli\`eres}, nous avons la
proposition suivante dont 
la preuve peut \^etre utile
dans d'autres contextes comme l'\'etude de la dimension de $\mu$ par exemple.
\begin{proposition} Soit $F$ une correspondance comme au
th\'eor\`eme 2.1. Alors il existe un ensemble pluripolaire $\E'$
tel que pour tout $z\not\in\E'$, tout $\epsilon>0$ et tout
$\delta>0$, la boule
$\Ball(z,r)$ poss\`ede au moins $(1-\epsilon)d_2^n$ branches
inverses r\'eguli\`eres
d'ordre $n$, de taille $(l^{-(1-\delta)m/2})$
o\`u $r=r(z,\delta,\epsilon)>0$ est une constante ind\'ependante de $n$.
\end{proposition}
\par
L'id\'ee de la d\'emonstration consiste \`a construire les
branches inverses r\'eguli\`eres
pour des unions de disques holomorphes
centr\'ees en un point g\'en\'erique $z$. Ensuite, utilisant un
th\'eor\`eme d'analyse complexe (lemme 3.8),
on peut prolonger les applications
holomorphes $F^{-n}_\B$
associ\'ees \`a ces branches inverses r\'eguli\`eres. Ces
applications sont d\'efinies aux voisinages de grandes familles de
disques centr\'es en $z$. On les prolonge en applications
holomorphes sur une petite boule $\Ball(z,r)$. Ces prolongements
fournissent les branches inverses r\'eguli\`eres pour $\Ball(z,r)$.
\par
La construction est faite par r\'ecurrence. Nous r\'esumons ici
le passage du rang $n-1$ au rang $n$. Soit $K_{-n+1}$ un disque
holomorphe fabriqu\'e au rang $n-1$
\`a partir d'un disque $K_0$ de rayon $r'>0$ centr\'e en $z$.
Nous voulons en fait construire des branches inverses
r\'eguli\`eres d'ordre 1
de $K_{-n+1}$.  D'abord, pour construire $\widehat K_{-n}$, on veut
que $K_{-n+1}$ ne rencontre pas les valeurs critiques de
$\pi_{2|Y}$; ensuite, pour obtenir $K_{-n}$, nous avons besoin que
$\pi_{1|\widehat K_{-n}}$ soit injective.
\par
Pour que ces deux conditions soient satisfaites, nous
s\'electionnons seulement les disques $K_{-n+1}$ de diam\`etre
assez petit qui ne sont pas trop proches d'un certain
sous-ensemble alg\'ebrique $(P=0)$ que nous appelons {\it
l'ensemble des valeurs critiques} de $F$. Pour chaque disque
s\'electionn\'e, on peut
construire $d_2$ branches inverses r\'eguli\`eres
d'ordre 1. Afin de continuer
la construction, nous devons s\'electionner, parmi les nouveaux
disques, ceux de petit diam\`etre qui ne sont pas trop proches de
$(P=0)$. L'hypoth\`ese sur l'exposant de Lojasiewicz implique
qu'un bon nombre de disques sont de petite aire. Ces disques ne
sont pas forc\'ement de petit diam\`etre. Mais, en diminuant
l\'eg\`erement $r'$, on rend leur diam\`etre petit.
Le fait que $\mu$ soit une mesure PB implique que pour $z$
g\'en\'erique, il n'y a qu'un petit
nombre de disques qui sont proches de $(P=0)$. 
\par
Soit $Q$ un polyn\^ome non nul.
Pour tout $\alpha>0$, notons
$\V_Q(\alpha)$ l'ensemble des
points $z$ v\'erifiant $|Q(z)|\leq\alpha$. Le lemme suivant
montre que les pr\'eimages de $F$ ne sont pas ``trop proches'' de
l'hypersurface $(Q=0)$. On peut le g\'en\'eraliser 
aux ensembles pluripolaires de 
$\P^k$. 
\begin{lemme} Soit $a>1$ une constante. Alors
il existe un ensemble pluripolaire $\E^a_Q$ tel que
$\sum_{n\geq 0}\mu^z_n(\V_Q(a^{-n}))$ soit fini pour
tout $z\not\in \E^a_Q$. 
\end{lemme}
\begin{preuve} On peut supposer $a<l^{1/2}$. 
On choisit une fonction $\chi$ lisse, positive sur
$\R^+$ telle que $\chi\leq 1$, $\chi=1$ sur $[0,1]$
et $\chi=0$ sur $[2,+\infty[$.
Soit $\tau\leq 1$ une fonction lisse, positive \`a support dans $V$ et
\'egale \`a 1 au voisinage de $\overline U$.
Posons $\psi^n:=\tau\chi(a^nQ)$.
Les mesures $\mu^z_n$ \'etant port\'ees par $\overline U$ pour
$n$ assez grand, il suffit de montrer que la
somme $\sum_{n\geq 0}\langle\mu^z_n,\psi^n\rangle$
est finie pour tout 
$z\not\in \E^a_Q$ o\`u l'ensemble 
$\E^a_Q$ sera pr\'ecis\'e dans la suite.
\par
D'apr\`es la proposition 2.7, la fonction $\log|Q|$ est
$\mu$-int\'egrable. On en d\'eduit que la somme $\sum_{n\geq
0}\langle\mu,\psi^n\rangle$ est finie. 
En effet, on a
\begin{eqnarray*}
\sum_{n\geq 0} \langle\mu,\psi^n\rangle
& \leq &  \sum_{n\geq 0} \mu(\V_Q(2a^{-n}))
=\int\Big(\sum_{n\geq 0} \ind_{\V_Q(2a^{-n})} \Big) \d\mu \\
& \leq & \int\frac{1}{\log a} \Big(\big|\log|Q|\big|+\log 2 +\log a
\Big)\d \mu.
\end{eqnarray*}
Posons $c_n:=\langle\mu,\psi^n\rangle$.
\par
Par d\'efinition de $\psi^n$, il existe une constante $c>0$ telle que
$\|\psi^n\|_{{\cal C}^2}\leq ca^{2n}$. Posons
$\varphi:=\log(1+\|z\|^2)$ et
$\psi^{n-}:=ca^{2n}\varphi-\psi^n$. Ce sont des fonctions p.s.h. \`a
croissance logarithmique.
Posons $c_\varphi:=\langle\mu,\varphi\rangle$ et
$c_n^-:=\langle\mu,\psi^{n-}\rangle$.
On a $c_n=ca^{2n}c_\varphi - c_n^-$.
Rappelons qu'on a suppos\'e $a<l^{1/2}$.
D'apr\`es
la proposition 2.7 et le corollaire 2.8, la s\'erie de fonctions
$$\Phi^+:=\sum_{n\geq 0} ca^{2n}(\Lambda^n\varphi -c_\varphi)$$
converge ponctuellement vers une fonction p.s.h. 
En particulier, $\Phi^+$ est 
localement born\'ee sup\'erieurement.
La proposition 2.7 et le corollaire 2.8 appliqu\'es aux
fonctions $\psi^{n-}$ impliquent aussi que la s\'erie de fonctions
$$\Phi^-:=\sum_{n\geq 0} (\Lambda^n\psi^{n-} -c_n^-)$$
converge ponctuellement vers une fonction p.s.h.
\par
Notons
$\E^a_Q$ l'ensemble o\`u $\Phi^-$ vaut $-\infty$.
C'est un ensemble pluripolaire. Pour $z\not\in\E^a_Q$, on a
\begin{eqnarray*}
\lefteqn{\sum_{n\geq 0} \langle\mu^z_n,\psi^n\rangle -
\sum_{n\geq 0} \langle\mu,\psi^n\rangle=}\\
& = & \sum_{n\geq 0} (\langle\mu^z_n,\psi^n\rangle-c_n)\\
& = & \sum_{n\geq 0} ca^{2n}(\langle\mu^z_n,\varphi\rangle-c_\varphi)
- \sum_{n\geq 0} (\langle\mu^z_n,\psi^{n-}\rangle-c_n^-)\\
& = & \Phi^+(z)-\Phi^-(z)<+\infty.
\end{eqnarray*}
Or la somme
$\sum \langle\mu,\psi^n\rangle$ est finie
comme on l'a montr\'e ci-dessus. La somme
$\sum \langle\mu^z_n,\psi^n\rangle$ est donc aussi finie.
\end{preuve}
\begin{lemme}
Il existe un polyn\^ome $P$ (ind\'ependant de $U$, $V$) et
une constante $A_1>0$ (d\'ependante de $U$, $V$) tels que pour tout
point $z\in V\setminus (P=0)$ la boule $\Ball_0:=\Ball(z,r_0)$
admette exactement $d_t$ branches inverses r\'eguli\`eres d'ordre
$1$ pour $F$, o\`u $r_0:=A_1|P(z)|$.
\end{lemme}
\begin{preuve}
Notons $Z:=|Y|$ le support de $Y$ et $m_2$
le degr\'e de l'application $\pi_{2|Z}$.
Soit $\Sigma_2$ l'ensemble des points $z\in \C^k$ tels que
$(\pi_{2|Z})^{-1}(z)$ contienne moins de $m_2$ points. C'est un
sous-ensemble alg\'ebrique de $\C^k$. Il existe donc un polyn\^ome
$P_1$ tel que $\Sigma_2\subset (P_1=0)$.
D\'efinissons de la m\^eme mani\`ere 
l'ensemble des valeurs critiques $\Sigma_1$ de $\pi_{1|Z}$ et posons
$\Sigma':=F(\Sigma_1)$. Soit $P_2$ un polyn\^ome tel que
$\Sigma'\subset (P_2=0)$. Posons $P=P_1^mP_2^m$ o\`u $m\geq 1$ est un
entier assez grand.
\par
Soit $A_1>0$ une constante assez petite.
Puisque $A_1$ est petite et $m$ est grand,
$(\pi_{2|Z})^{-1}\Ball_0$ contient exactement $m_2$ composantes
connexes et la restriction de $\pi_1$ \`a chacune de ces
composante est injective.
\par
Soit $\widehat{\Ball}_{-1}$
une composante connexe de $\pi_2^{-1}(\Ball_0)\cap Y_i^*$,
pour une certaine indice $i$.
La suite
$$\Ball_{-1},(\widehat{\Ball}_{-1},i_1),\Ball_0$$
avec $i_1:=i$ et $\Ball_{-1}:=\pi_1(\widehat{\Ball}_{-1})$,
est une branche inverse r\'eguli\`ere d'ordre 1
de $\Ball_0$. Il y en a exactement $d_2$
branches.
\end{preuve}

Notons $\E_1$ l'orbite de $(P=0)$ par $(F^n)_{n\geq 0}$.
Par d\'efinition, tout point $z\not\in\E_1$ admet $d_2^n$
branches inverses d'ordre $n$.
Posons
$\E':=\E^a_P\cup \E_1$. Il est clair que $\E'$ est pluripolaire.
\par
Fixons maintenant des constantes $\delta$ avec $0<\delta<1$, 
$a$ avec $1<a<l^{(1-\delta)/2}$, 
$\epsilon$ avec $0<\epsilon<1$ et un point $z\in
V\setminus\E'$. Nous allons montrer
que la boule $\Ball(z,r)$ poss\`ede
au moins $(1-\epsilon)d_2^n$ branches inverses r\'eguli\`eres
d'ordre $n$ et de taille
$(l^{-(1-\delta)m/2})$ pour $r>0$ assez petit.
\par
Fixons
une droite $\Delta$ passant par $z$. Notons $\Delta_R$ le
disque de centre $z$ et de rayon $R$ dans $\Delta$. 
Rappelons le lemme de comparaison aire-diam\`etre d\^u \`a
Briend-Duval \cite{BriendDuval2}. Ce lemme est valable pour un
cas plus g\'en\'eral. Dans le cas pr\'esent, on peut le montrer
en utilisant la formule de Cauchy.
\begin{lemme} Soient $\pi$ une application holomorphe du disque
unit\'e $\Disc:=\Disc(0,1)$
dans $U$ et $\tau$ une constante v\'erifiant
$0<\tau<1$. Alors il existe une constante $A_2>0$ ind\'ependante de
$\pi$ et de $\tau$ telle que le diam\`etre de $\pi(\Disc(0,1-\tau))$
soit plus petit ou \'egal \`a
$A_2\sqrt{\tau^{-1}\aire(\pi(\Disc))}$, les points
de $\pi(\Disc)$ \'etant compt\'es avec multiplicit\'e.   
\end{lemme}
\par
Nous allons montrer la proposition suivante dans laquelle la
constante $A_3>1$ sera donn\'ee dans le lemme 3.8.
\begin{proposition} Pour tout $\epsilon_1>0$, il
existe $r>0$ et $n_0\geq 0$ ind\'ependants de $\Delta$ 
tels que $\Delta_{A_3r}$ poss\`ede au moins
$(1-\epsilon_1)d_t^n$ branches inverses r\'eguli\`eres
$\B$ d'ordre $n$,
de taille $(\frac{1}{2}l^{-(1-\delta)m/2})$ et telles que
$|P(F^{-m}_\B(z))|\geq a^{-m}$ pour tout $n_0\leq m\leq n$.
\end{proposition}
\par
Fixons $n_0\geq 1$ assez grand tel que les propri\'et\'es suivantes
soient satisfaites
\begin{enumerate}
\item $A_1a^{-n}> l^{-(1-\delta)n/2}$ pour tout $n\geq n_0$.
\item $\sum_{n\geq n_0+1} \mu_n^z(\V_P(a^{-n}))\leq \frac{1}{2}
\epsilon_1$.
\item $\sum_{n\geq n_0+1} n^2 l^{-\delta n}\leq
\frac{1}{8}\epsilon_1A_2^{-2}$.
\end{enumerate}
Posons $\nu_m:=\sum_{n=n_0+1}^m \mu^z_n(\V_P(a^{-n}))$ et
$\delta_m:=4A_2^2\sum_{n=n_0+1}^m n^2 l^{-\delta n}$.
Fixons $r_1>0$ assez
petit tel que pour tout $n\leq n_0$,
$\Ball(z,r_1)$ poss\`ede exactement $d_2^n$
branches inverses r\'eguli\`eres d'ordre $n$ de taille 
$(\frac{1}{2}l^{-(1-\delta)m/2})$.
Posons pour tout $n\geq n_0+1$
$$r_n:=\prod_{s=n_0+1}^n\left(1-\frac{1}{s^2}\right) r_1.$$
Cette suite $(r_n)$ d\'ecro\^{\i}t vers une constante $A_3r>0$.
Il suffit pour la proposition
3.6 de montrer par r\'ecurrence sur $n\geq n_0$
que $\Delta_{r_n}$ poss\`ede au
moins $(1-\nu_n-\delta_n) d_2^n$ branches inverses r\'eguli\`eres
$\B$ d'ordre $n$,
de taille $(\frac{1}{2}l^{-(1-\delta)m/2})$ et telles que
$|P(F^{-m}_\B(z))|\geq a^{-m}$ pour tout $n_0\leq m\leq n$.
Supposons que c'est vrai au rang
$n-1\geq n_0$.
Montrons le au rang $n$. Notons $\F$ la famille des
branches inverses r\'eguli\`eres
d'ordre $n-1$  de $\Delta_{r_{n-1}}$
v\'erifiant la proposition 3.6.
\par
Si $\B$ est un \'el\'ement de $\F$, d'une part, le point
$w:=F^{-n+1}_\B(z)$ v\'erifie $|P(w)|\geq a^{-n+1}$, d'autre part,
l'ensemble $W:=F^{-n+1}_\B(\Delta_{r_{n-1}})$
est contenu dans la boule
$\Ball(w,A_1a^{-n+1})$ car son diam\`etre est plus petit que
$\frac{1}{2}A_1a^{-n+1}$. Or, d'apr\`es le lemme 3.4, 
cette boule $\Ball(w,A_1a^{-n+1})$ admet exactement
$d_2$ branches inverses r\'eguli\`eres d'ordre 1. On en d\'eduit que
$W$ admet aussi $d_2$ branches inverses r\'eguli\`eres
d'ordre 1. Ceci est vrai
pour tout $\B\in \F$. En somme, $\Delta_{r_{n-1}}$
poss\`ede au moins $(1-\nu_{n-1}-\delta_{n-1})d_2^n$
branches inverses r\'eguli\`eres
d'ordre $n$. Notons $\G$ cette famille de branches inverses
r\'eguli\`eres d'ordre $n$. On a le lemme suivant.
\begin{lemme} La somme $\sum \aire(F^{-n}_\B(\Delta_{r_{n-1}}))$ 
pour $\B\in\G$ est plus petite que
$d_2^nl^{-n}$. 
\end{lemme}
\begin{preuve} Notons $[\Delta]$ le courant d'int\'egration sur
$\Delta$. Il suffit de majorer la masse de $(F^n)^*[\Delta]$
par $d_2^nl^{-n}$. On a d'apr\`es un lemme de comparaison
(voir \cite{Sibony} ou \cite[proposition 5.4]{DinhSibony2}) 
\begin{eqnarray*}
\|(F^n)^*[\Delta]\| & = & \langle (F^n)^*[\Delta],\omega\rangle
= \langle[\Delta], (F^n)_*\omega\rangle \\
& = & \frac{1}{2}d_2^n l^{-n}\int_\Delta \ddc
(l^n\Lambda^n \log(1+\|z\|^2)) \\
& \leq & \frac{1}{2}d_2^n l^{-n} \int_\Delta \ddc\log(1+\|z\|^2)\\
& = & d_2^n l^{-n} \int_\Delta \omega
=d_2^nl^{-n}
\end{eqnarray*}
L'in\'egalit\'e ci-dessus est une cons\'equence du fait que
$l^n\Lambda^n\log(1+\|z\|^2)-\log(1+\|z\|^2)$
est born\'ee sup\'erieurement.
\end{preuve}
{\bf Fin de la d\'emonstration de la proposition 3.6.}
Notons $\G_1$ l'ensemble des branches inverses r\'eguli\`eres
$\B\in\G$ telles
que l'aire de $F^{-n}_\B(\Delta_{r_{n-1}})$ exc\`ede
$\frac{1}{4}A_2^{-2}n^{-2}l^{-(1-\delta)n}$. 
D'apr\`es le lemme 3.7, $\G_1$ contient au plus
$4A_2^2n^2l^{-\delta n}d_2^n$ \'el\'ements.
Par cons\'equent, le cardinal de la famille $\G':=\G\setminus \G_1$
est au moins \'egal \`a
$$\#\G -\# \G_1\geq (1-\nu_{n-1}-\delta_{n-1})d_2^n -
4A_2^2n^2l^{-\delta n}d_2^n = (1-\nu_{n-1}-\delta_n)d_2^n.$$
D'apr\`es le lemme 3.5 (appliqu\'e \`a l'application $F^{-n}_\B$
sur les disques $\Delta_{r_{n-1}}$ et $\Delta_{r_n}$), pour
tout $\B\in\G'$, le diam\`etre de
$F^{-n}_\B(\Delta_{r_n})$ est au plus \'egal \`a
$\frac{1}{2}l^{-(1-\delta)n/2}$.
\par
Notons $\G_2$ la famille des \'el\'ements $\B\in \G'$ tels que
$F^{-n}_\B(z)\in \V_P(a^{-n})$. Alors $\G_2$ contient au plus
$\mu^z_n(\V_P(a^{-n}))d_2^n$ \'el\'ements. Par
cons\'equent, la famille 
$\G'':=\G'\setminus\G_2$ contient au moins
$$ (1-\nu_{n-1}-\delta_n)d_2^n - 
\mu^z_n(\V_P(a^{-n})) d_2^n
=  (1-\nu_n-\delta_n) d_2^n$$
\'el\'ements.
Cette famille $\G''$ de branches
inverses r\'eguli\`eres d'ordre $n$ de $F$ v\'erifie
la proposition 3.6.  
\par
\hfill $\square$
\\
{\bf Fin de la d\'emonstration de la proposition 3.2.}
Prenons $\epsilon_1=\epsilon/2$ et 
notons $\W$ la famille des 
droites complexes passant par $z$. Cette famille
est param\'etr\'ee par 
l'espace projectif $\P^{k-1}$. Notons $\H_{2k-2}$
la mesure de Hausdorff $(2k-2)$-dimensionnelle de 
masse 1 sur $\W$. Notons \'egalement
$\F_n$ la famille des branches
inverses r\'eguli\`eres d'ordre $n$ de $z$. 
Cette famille contient exactement $d_2^n$ \'el\'ements
car $z\not\in \E_1$. 
Pour tout
$\B_z\in \F_n$ notons $\W_{\B_z}$ la famille 
des droites $\Delta\in \W$ telles que $\Delta_{A_3r}$
poss\`ede une branche inverse r\'eguli\`ere d'ordre $n$,
accroch\'ee \`a la branche $\B_z$ et v\'erifiant les
propri\'et\'es dans la
proposition 3.6. D'apr\`es cette proposition,
on a $\sum_{\B_z} \H_{2k-2}(\W_{\B_z})\geq 
(1-\epsilon_1) d_2^n$ lorsque $r>0$ est suffisamment petit.
Notons $\F_n'$ la famille des $\B_z$ tels que 
$\H_{2k-2}(\W_{\B_z})\geq 1/2$.
Du fait que $\H_{2k-2}(\W_{\B_z})\leq 1$ 
pour tout $\B_z$, on a
$$\# \F_n' + \frac{1}{2}\big(d_2^n-\#\F'_n\big) 
\geq \sum \H_{2k-2}(\W_{\B_z})\geq 
(1-\epsilon_1) d_2^n.$$
On en d\'eduit que $\#\F_n'\geq (1-2\epsilon_1)d_2^n=
(1-\epsilon)d_2^n$.
On va appliquer le th\'eor\`eme de Sibony-Wong 
suivant pour chaque $\W_{\B_z}$ avec  
$\B_z\in\F_n'$.
\begin{lemme}{\bf \cite{Alexander, SibonyWong}}
Soit $\alpha>0$ une constante positive.
Soit $A_3>1$ une constante suffisamment grande et soit
$\W'$ une famille de droites passant par $z$. 
Supposons que $\H_{2k-2}(\W')\geq \alpha$. Notons
$\Sigma$ l'intersection 
de ces droites avec
la boule $\Ball(z,A_3r)$. Alors toute application
holomorphe $f$ d'un
voisinage de $\Sigma$ d'image dans $\C^k$ se prolonge 
en application holomorphe de
$\Ball(z,r)$ dans $\C^k$. De plus, 
on a
$$\sup_{w\in \Ball(z,r)}\|f(w)-f(z)\| \leq
\sup_{w\in \Sigma} \|f(w)-f(z)\|.$$
En particulier, on a $\diam f(\Ball(z,r))\leq 2\diam f(\Sigma)$.
\end{lemme}
\par
On prend $\alpha=1/2$. 
Fixons un $\B_z\in\F'_n$. Notons $K_0$ l'intersection des droites de 
$\W_{\B_z}$ avec $\Ball(z,A_3r)$ et $\B$
sa branche inverse r\'eguli\`ere d'ordre $n$ 
accroch\'ee \`a 
$\B_z$. D'apr\`es le lemme 3.8, l'application $F^{-m}_{\B}$,
qui est holomorphe
au voisinage de $K_0$, se prolonge en
une application holomorphe de $\Ball(z,r)$
dans $\C^k$ pour tout
$1\leq m\leq n$. De plus, son image est de diam\`etre au plus 
$l^{-(1-\delta)m/2}$. 
Les applications obtenues sont injectives. En effet, on montre
par r\'ecurrence que l'image de $\Ball(z,r)$ par $F_\B^{-m}$ est
contenue dans la boule de rayon $A_1|P(F_\B^{-m}(z))|$ centr\'ee
en $F_\B^{-m}(z)$ (voir le lemme 3.4). Chaque prolongement
holomorphe
fournit une branche inverse r\'eguli\`ere d'ordre 
$n$, de taille $(l^{-(1-\delta)m/2})$ 
pour la boule $\Ball(z,r)$. Ceci termine la preuve
de la proposition 3.2.
\par
\hfill $\square$
\par

La suite de ce paragraphe a \'et\'e d\'emontr\'ee en collaboration avec
Charles Favre.
Notons $\PC_1$ l'ensemble des points $z\in\C^k$ tels que au moins un des 
germes locaux irr\'eductibles de $Y$ en 
$\pi_2^{-1}(z)\cap Y$ ne se projete pas injectivement par $\pi_2$ sur $\C^k$. 
C'est une hypersurface de $\C^k$.
Soit
$S_1$ le courant d'int\'egration sur $\PC_1$. Posons
$S_n:=(F^{n-1})_*(S_1)$ et $S:=\sum_{n\geq 1}d_2^{-n+1}S_n$. 
Ce sont des courants positifs ferm\'es de bidegr\'e $(1,1)$ de
$\C^k$. Le courant $S$ est bien d\'efini car 
d'apr\`es le lemme 2.2, la masse de $d_2^{-n} S_n$ est de l'ordre de
$l^{-n}$. Soient
$\PC_n:=\cup_{i=1}^{n-1} \supp(S_i)$ et $\PC_\infty:=\supp(S)$
{\it les ensembles postcritiques d'ordre $n$ et d'ordre infini}
de $F$.
Posons \'egalement $\Sigma:=Y\cap
\pi_2^{-1}(\PC_1)$.
\par
Soit $K_0$ un disque holomorphe de centre $z$,
une famille de disques holomorphes centr\'es
en $z$ ou une boule de centre $z$. On ne consid\`ere
que les disques holomorphes plats
qui ne sont pas contenus dans $\PC_\infty$. On appelle {\it
branche inverse d'ordre $n$ de $K_0$} (voir Figure 2) toute suite $\B$
$$K_{-n}, (\widehat K_{-n},i_n), K_{-n+1}, \ldots, K_{-1}, (\widehat
K_{-1},i_1), K_0$$
munie des applications holomorphes $\widehat
F_\B^{-m}:K_0\longrightarrow \widehat K_{-m}$
telle que
\begin{enumerate}
\item[(i)]
Les ensembles 
$K_{-m}\subset\C^k$, $\widehat K_{-m}\subset
Y_{i_m}^*\setminus \Sigma$ sont connexes;
$\pi_1(\widehat K_{-m})=K_{-m}$,
$\pi_2(\widehat K_{-m})=K_{-m+1}$ pour $1\leq m\leq n$;
\item[(ii)]
$\pi_2\circ\widehat F_\B^{-1}=\id$ et
$\pi_1\circ \widehat F^{-m+1}_\B = \pi_2\circ \widehat
F^{-m}_\B$ pour tout $2\leq m\leq n$.
\end{enumerate}
On n'exige pas que $\pi_{1|\widehat K_{-n}}$ soit injective.
Posons $F^{-m}_\B:=\pi_1\circ\widehat F^{-m}_\B$.
\begin{figure}
\begin{center}
\psfrag{A}{$\C^k$}
\psfrag{B}{$K_0$}
\psfrag{C}{$K_{-1}$}
\psfrag{D}{$K_{-2}$}
\psfrag{E}{$\widehat K_{-1}$}
\psfrag{F}{$\widehat K_{-2}$}
\includegraphics{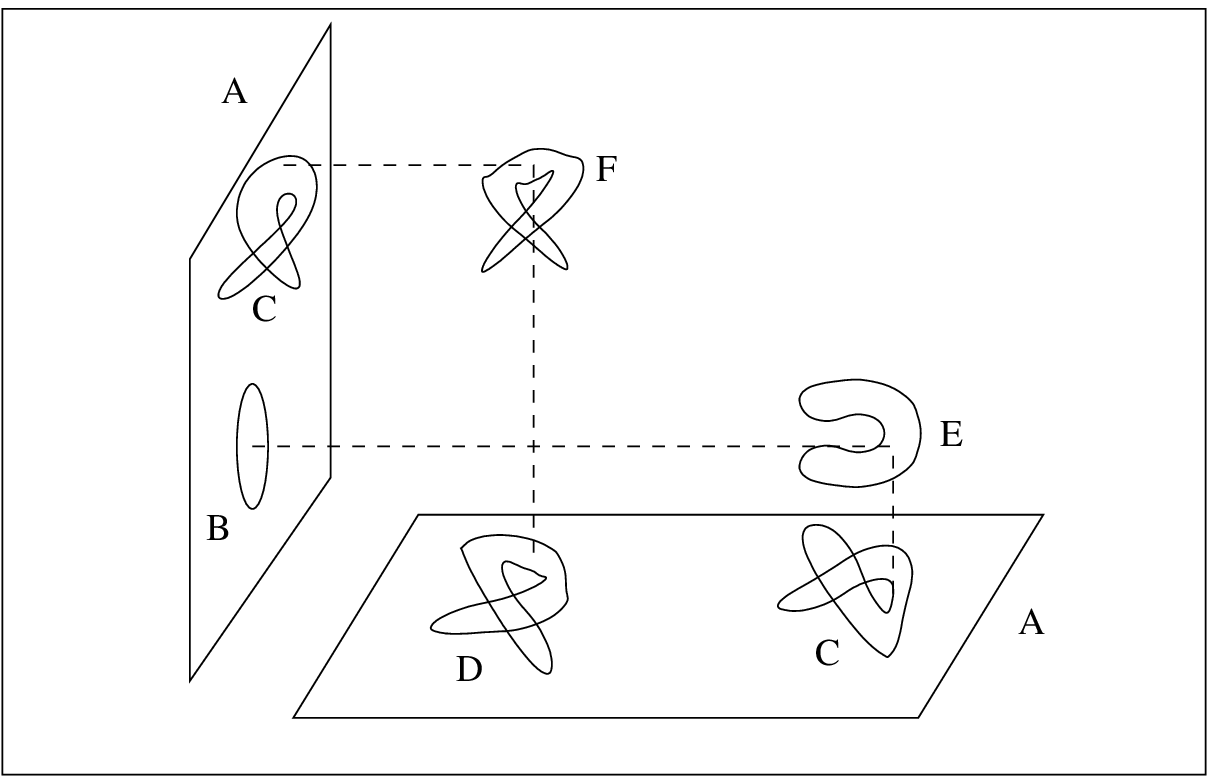}
\caption{Branche inverse non r\'eguli\`ere.}
\end{center}
\end{figure}
\par
Fixons un point $w\in K_0$.
La branche inverse $\B$ de $K_0$
est accroch\'ee \`a la branche inverse r\'eguli\`ere $\B_w$ de
$w$ donn\'ee par la suite 
$$w_{-n},(\widehat w_{-n},i_n), w_{-n+1}, (\widehat
w_{-n+1},i_{n-1}), \ldots,
w_{-1}, (\widehat w_{-1},i_1), w_0$$
o\`u $w_0:=w$, $\widehat w_{-m}:=\widehat F_\B^{-m}(w_0)$ et
$w_{-m}:=\pi_1(\widehat w_{-m})$. Par unicit\'e du prolongement
analytique, la branche $\B$ est uniquement d\'etermin\'ee par la
branche $\B_w$. Autrement dit, si deux branches inverses d'ordre
$n$ de $K_0$
sont accroch\'ees \`a une m\^eme branche r\'eguli\`ere
$\B_w$, alors elles
sont \'egales.
On en d\'eduit que $K_0$ poss\`ede au plus
$d_2^n$ branches inverses d'ordre $n$. Nous avons le th\'eor\`eme
suivant.
\begin{theoreme} Soit $F$ une correspondance comme au
th\'eor\`eme 2.1. Soit $z\in\C^k$ un point tel que
le nombre de Lelong de $S$ en $z$ v\'erifie
$\nu(S,z)<1$. Alors pour tout $\epsilon>0$ et tout $\delta>0$, la
boule $\Ball(z,r)$ poss\`ede au moins
$(1-\nu(S,z)-\epsilon)d_2^n$ branches inverses d'ordre $n$ de
taille $(l^{-(1-\delta)m/2})$ o\`u $r>0$ est une constante
ind\'ependante de $n$.  
\end{theoreme}
Notons que ce th\'eor\`eme a \'et\'e prouv\'e dans \cite{DinhSibony1} pour les 
endomorphismes holomorphes de degr\'e $d\geq 2$ 
de $\P^k$ (ou plus g\'en\'eralement pour certaines 
applications d'allure polynomiale). Dans ce cas, les branches inverses sont de taille
$\simeq (d^{-m/2})$. Ind\'ependemment, dans une note non publi\'ee, 
Briend-Duval ont construit les branches inverses sur les boules
par une autre m\'ethode. Ils ont montr\'e que la taille de ces branches tend
vers 0 quand $n\rightarrow\infty$. Les branches inverses sur les disques ont \'et\'e
d\'ej\`a construites dans \cite{BriendDuval2}. Le lecteur trouvera dans \cite{Dinh2}
d'autre version du th\'eor\`eme 3.9 avec des applications.

Pour prouver ce th\'eor\`eme,
il suffit de construire des branches inverses
pour des disques holomorphes plats centr\'es en $z$.  Les deux
arguments suivants permettent d'adapter la preuve de la proposition
3.2.
\par
1. Fixons $r'>0$ assez petit. Notons $s_n$ le nombre de points
d'intersection de $\Delta_{r'}$ avec $F^{n-1}(\PC_1)$, \cad la masse
de la mesure d'intersection de $[\Delta_{r'}]$ avec
le courant $S_n$.
Ce nombre $s_n$ est aussi \'egal \`a la masse de la mesure
d'intersection
de $S_1$ avec $(F^{n-1})^*[\Delta_{r'}]$. 
Notons $\W'$
la famille des droites $\Delta$ telles que $\sum d_2^{-n+1} s_n
\leq \nu(S,z)+\epsilon$. Par tranchage, si $r'$ est suffisamment
petit, on a $\H_{2k-2}(\W')\geq 1/2$.
\par
2. Fixons un $\Delta$ dans $\W'$. On construit par r\'ecurrence
sur $n$, au moins $N_n:=d_2^n (1-\sum_{m=1}^n d_2^{-m+1}
s_m)$ branches inverses d'ordre $n$ pour $\Delta_{r'}$. R\'esumons
le passage du rang $n-1$ au rang $n$. Par d\'efinition de $s_n$,
il existe au moins $N_{n-1}- s_n$
branches inverses $\B$ d'ordre $n-1$ de $\Delta_{r'}$
$$K_{-n+1},(\widehat K_{-n+1},i_{n-1}),
\ldots, K_{-1}, (\widehat K_{-1},i_1), K_0$$
avec $K_0=\Delta_{r'}$
telles que
l'ensemble
$F^{-n+1}_\B(\Delta_{r'})$, qui est \'egal \`a $K_{-n+1}$,
ne rencontre pas $\PC_1$.
D'autre part, l'application $\pi_2$ d\'efinit un rev\^etement
non ramifi\'e de degr\'e $d_2$ de
$Y\setminus \Sigma$ au dessus de
$\C^k\setminus \PC_1$. L'ensemble $\Delta_{r'}$ \'etant simplement 
connexe, on peut construire, pour une telle branche $\B$,
exactement $d_2$ applications
holomorphes 
$\tau$ de $\Delta_{r'}$ dans un des $Y^*_i$
telles que $\pi_2\circ\tau = F^{-n+1}_\B$.
On obtient donc au moins $N_n$ branches inverses $\B'$ d'ordre
$n$ de $\Delta_{r'}$
$$K_{-n},(\widehat K_{-n},i_n),\ldots, K_{-1},
(\widehat K_{-1},i_1), K_0$$
en posant $i_n:=i$, $\widehat F^{-n}_{\B'}:=\tau$, $\widehat
K_{-n}:=\tau(\Delta_{r'})$ et $K_{-n}:=\pi_1(\widehat K_{-n})$.
\par
\ \hfill $\square$
\begin{corollaire} Soit $F$ une correspondance comme au
th\'eor\`eme 2.1. Alors l'ensemble exceptionnel $\E$ de $F$ est
\'egal \`a $\cup_{n\geq 0} F^n(\E_0)$ o\`u $\E_0$
est le plus grand sous-ensemble alg\'ebrique propre de $\C^k$
v\'erifiant $F^{-1}(\E_0)=\E_0$.
\end{corollaire}
\begin{preuve}
Pour tout $\nu>0$,
notons $X_\nu$ l'ensemble des points $z$ tels que $\nu(S,z)\geq
\nu$. C'est un sous-ensemble alg\'ebrique de $\C^k$ contenu dans
$\PC_\infty$. Soient $z\not\in X_\nu$ avec $\nu<1$, $\epsilon$ et
$\Ball(z,r)$ v\'erifiant le th\'eor\`eme 3.9.
D'apr\`es la proposition 3.1,
pour un point g\'en\'erique $w\in \Ball(z,r)$, on a $\lim
\mu^w_n=\mu$. 
Consid\'erons une valeur adh\'erente $\mu^z$
de la suite de mesures $(\mu^z_n)$ et $\mu^z_\reg$ sa partie
absolument continue par rapport \`a $\mu$.
D'apr\`es le th\'eor\`eme 3.9, $z$ et $w$ poss\`edent au moins
$d_2^n(1-\nu-\epsilon)$ images r\'eciproques $z_{-n}$ et
$w_{-n}$
d'ordre $n$ telles que la distance entre $z_{-n}$ et $w_{-n}$
soit major\'ee par $l^{-(1-\delta)n/2}$.
On en d\'eduit que la masse de $\mu^z_\reg$ est au moins \'egale
\`a $1-\nu-\epsilon$.
Donc elle est au moins \'egale \`a $1-\nu$.
On d\'eduit aussi que $\E\subset \cup_{\nu>0} X_\nu \subset\PC_\infty$. 
\par
Notons $\E_{X_\nu}$ l'ensemble des points $z$ tels que
$F^{-n}(z)\subset X_\nu$ pour tout $n\geq 0$.
C'est le plus grand sous-ensemble alg\'ebrique de
$X_\nu$ qui v\'erifie $F^{-1}(\E_{X_\nu})\subset
\E_{X_\nu}$. Il est clair que $\E_{X_\nu}\subset \E$. 
Pour tout
$z\in \C^k$, posons $\tilde\mu^z_0:=\ind_{X_\nu}\delta_z$,
$\tilde\mu^z_n:=d_2^{-1}\ind_{X_\nu} F^*(\tilde\mu^z_{n-1})$.
Les masses de $\tilde\mu_n^z$ d\'ecroissent vers une constante
$\tau_\nu(z)$. De la m\^eme mani\`ere que dans
\cite[3.4.1-3.4.3]{DinhSibony1}, on montre que
$\tilde\E_{X_\nu}:=\cup_{n\geq 0}F^n(\E_{X_\nu})$ contient
l'ensemble 
$(\tau_\nu>0)$.
La mesure $\mu^z_n$ \'etant \'egale \`a
la moyenne des mesures $\mu^a_m$ avec
$a\in F^{-n+m}(z)$ et $m\leq n$,
si $\tau_\nu(z)=0$ la masse de $\mu^z_\reg$ est, comme pour
$z\not\in X_\nu$,
au moins \'egale \`a $1-\nu$. En particulier,
si $z$ n'appartient pas \`a
$\cup_{\nu>0} \tilde\E_{X_\nu}$, $\mu^z_\reg$ est de masse totale.
Dans ce cas,
la mesure $\mu^z$ est \'egale \`a $\mu^z_\reg$ et ne charge
pas $\PC_\infty$. Ceci implique que $\mu^z_n(\PC_\infty)$ tend
vers $0$. Or $\E\subset\PC_\infty$. On d\'eduit qu'alors
$\lim\mu^z_n=\mu$ et par cons\'equent $\E=\cup_{\nu>0} \tilde
\E_\nu$. 
\par
Montrons que $\E_{X_\nu}=\E_{X_1}$ pour tout $0<\nu<1$. Il est
clair que $\E_{X_\nu}\supset\E_{X_1}$. Observons
que d'apr\`es le th\'eor\`eme 3.9,
si $z\not\in \E_{X_1}$ on a $\mu^z_\reg\not=0$. Si $z$
appartient \`a $\E_{X_\nu}$, les mesures $\mu^z_n$ sont port\'ees
par $X_\nu$. Donc $\mu^z_\reg=0$ car $\mu$ ne charge
pas les ensembles analytiques et $\mu^z_n$ est absolument
continue par rapport \`a $\mu$. D'o\`u on d\'eduit que $z\in
\E_{X_1}$.
\par
On a montr\'e que $\E=\tilde \E_{X_1}$. Posons $\E_0:=\cap_{n\geq
0}F^{-n}(\E_{X_1})$. C'est un sous-ensemble alg\'ebrique de $\C^k$
v\'erifiant $F^{-1}(\E_0)=\E_0$. De plus,
$\E_0$ est \'egal \`a
l'intersection d'une famille finie d'ensembles
alg\'ebriques $F^{-n}(\E_{X_1})$.
On a donc $\E=\tilde
\E_{X_1}=\cup_{n\geq 0}F^n(\E_0)$. Si $\E_0'$ est un
sous-ensemble alg\'ebrique v\'erifiant $F^{-1}(\E_0')=\E_0'$, on
montre comme on l'a fait pour $\E_{X_\nu}$
que $\E_0'\subset \E_{X_1}$. La propri\'et\'e
$F^{-1}(\E_0')=\E_0'$ implique que $\E_0'\subset
\E_0=\cap_{n\geq 0} F^{-n}(\E_{X_1})$. D'o\`u on d\'eduit que
$\E_0$ est le plus grand sous-ensemble alg\'ebrique propre de
$\C^k$ invariant par $F^{-1}$.
\end{preuve}
\par
Lorsque $F$ est une application polynomiale \cite{DinhSibony1},
on a $\E=\E_0$; l'existence de points p\'eriodiques r\'epulsifs
implique que
la pseudo-m\'etrique de Kobayashi de $\C^k\setminus \E_0$ est
identiquement nulle; si $X$ est un sous-ensemble alg\'ebrique
v\'erifiant $F^{-1}(X)\subset X$ alors $F^{-1}(X)=X$ et $F(X)=X$.
Toutes ces
propri\'et\'es sont fausses en g\'en\'eral pour les
correspondances comme on le voit dans les exemples suivants. 
\begin{exemples} \rm Les correspondances suivantes
sont polynomiales sur $\C$. Leurs exposants de Lojasiewicz sont
strictement sup\'erieurs \`a 1.
\par
Consid\'erons
les polyn\^omes d'une variable
$f(z):=z^3$ et $g(z):=z^2-z$ et la correspondance 
$F:=f\circ g^{-1}$ sur $\C$ associ\'ee \`a la courbe
$Y=\{(g(z),f(z)),\ z\in\C\}$ de $\C^2$. On v\'erifie
que $\E_0=\{0\}$ et  
$\E=\cup_{n\geq 0} (f\circ g^{-1})^n(0)$. On v\'erifie aussi que
$0$ n'est pas isol\'e dans $\E$ et 
donc l'ensemble $\E$ est parfait.
\par
Notons $(z_1,z_2)$ les coordonn\'ees de $\C^2$. Soient $f$, $g$
deux polyn\^omes d'une variable tels que $1<\deg(g)<\deg(f)$.
Consid\'erons la
correspondance $F$ donn\'ee par la courbe
$Y:=\{[f(z_1)]^2+f(z_2)g(z_2)=0\}$.
On a
$F^{-1}(f=0)=(f=0)$ et donc $(f=0)\subset \E_0$.
L'ouvert
$\C\setminus \E_0$ est Kobayashi hyperbolique
lorsque $f$ poss\`ede deux racines distinctes. 
\end{exemples}
\begin{remarque} \rm
On peut \'etudier la restriction de la correspondance $F$ \`a
$\E_0$, \cad la correspondance sur $\E$ dont le graphe est \'egal
\`a $Y\cap (\E_0\times\E_0)$,
et construir sur $\E_0$ une mesure invariante
par $F^*$ qui ne charge pas les sous-ensembles
pluripolaires de $\E_0$. En faisant une r\'ecurrence
descendante sur la dimension et sur le nombre de composantes du
graphes, on montre que le c\^one des mesures positives $\nu$
qui v\'erifient $F^*(\nu)=d_2\nu$, est de dimension fini.
Pour tout $z\in\C^k$,
les valeurs adh\'erentes \`a la suite $(\mu^z_n)$
sont singuli\`eres par rapport \`a $\mu$ si et seulement si  
$z\in\E_0$ (voir aussi \cite{DinhSibony4, Dinh2}). 
\end{remarque}
\section{Points p\'eriodiques r\'epulsifs}
On appelle {\it point fixe} de $F$ tout point $z$ appartenant \`a 
l'ensemble $\pi_2(Y\cap \Diag)$ o\`u $\Diag$ est la diagonale de 
$\C^k\times \C^k$. Les points
fixes isol\'es de $F$ sont compt\'es avec
multiplicit\'es. {\it Les points p\'eriodiques de p\'eriode $n$} 
sont les points fixes de $F^n$. On dira que $z_0$ est {\it 
un point p\'eriodique 
r\'egulier de p\'eriode $n$} s'il existe une branche
inverse r\'eguli\`ere d'ordre $n$ d'un 
voisinage ouvert $K_0$ de $z_0$
$$K_{-n}, (\widehat K_{-n},i_n), K_{-n+1}, (\widehat
K_{-n+1},i_{n-1}), \ldots, 
K_{-1}, (\widehat K_{-1},i_1), K_0$$
qui est accroch\'ee \`a une branche inverse r\'eguli\`ere
$\B_0$ d'ordre $n$ de $z$
donn\'ee par la suite 
$$z_{-n}, (\widehat z_{-n},i_n), z_{-n+1}, (\widehat
z_{-n+1},i_{n-1}), 
\ldots, z_{-1}, (\widehat z_{-1},i_1), z_0$$
avec $z_{-n}=z_0$.
Si, de plus, les valeurs propres de
la d\'eriv\'ee de $F^{-n}_{\B_0}$ en $z_0$
sont de module strictement 
plus petit que $1$, on dit que $z_0$ est {\it 
p\'eriodique r\'egulier r\'epulsif}. La proposition
suivante donne le nombre 
de points p\'eriodiques (voir \cite{Dinh2} pour le cas des correspondances sur les
vari\'et\'es k\"ahl\'eriennes compactes).
\begin{proposition} Soit $F$ une correspondance
comme au th\'eor\`eme 2.1.
Alors $F$ admet
exactement $d_2^n$ points p\'eriodiques de p\'eriode $n$ 
compt\'es avec multiplicit\'es. Si $G$ est une autre
correspondance polynomiale
de degr\'e topologique $(p_1,p_2)$ alors pour $n$
assez grand $G\circ F^n$ admet exactement $p_2d_2^n$ points fixes
compt\'es avec multiplicit\'es. 
\end{proposition}
\begin{preuve} 
Puisque l'exposant de Lojasiewicz de $F$ est sup\'erieur \`a 1,
le choix de la boule $V$ entra\^{\i}ne que $(\Diag\cap Y)\subset
(V\times V)$.
De plus, on a $Y\cap (bV\times V)=\emptyset$. 
Pour tout $0\leq t\leq 1$, posons 
$$\Diag(t):=\{(x,tx) \mbox{ avec } x\in\C^k\}.$$
C'est un sous-espace de 
$\C^k\times \C^k$ et on a $\Diag(1)=\Diag$. Comme   
$Y\cap (bV\times V)=\emptyset$, le nombre de
points d'intersection de 
$Y\cap (V\times V)$ avec 
$\Diag(t)$ ne d\'epend pas de $t$. Quand $t=0$,
ce nombre est exactement 
le nombre de pr\'eimages de $0$ dans $V$ compt\'e
avec multiplicit\'es. Il 
est donc \'egal \`a $d_2$. Par cons\'equent, $F$ admet exactement
$d_2$ points fixes. De m\^eme, on montre que $F^n$ admet
exactement $d_2^n$
points fixes. Ce sont les points p\'eriodiques de
p\'eriode $n$ de $F$.
\par
Pour $n$ assez grand,
l'exposant de Lojasiewicz de $G\circ F^n$ est strictement plus
grand que $1$. Par cons\'equent,
$G\circ F^n$ admet exactement $p_2d_2^n$
points fixes. 
\end{preuve} 
\par
Le th\'eor\`eme suivant g\'en\'eralise un r\'esultat de
Lyubich \cite{Lyubich} qui
a consid\'er\'e le cas des fractions rationnelles de
$\P^1$.
\begin{theoreme} Soient $F$ une correspondance 
comme au th\'eor\`eme 2.1 et $G$ une autre
correspondance polynomiale de degr\'e topologique $(p_1,p_2)$ et
d'exposant de Lojasiewicz quelconque.
Notons $\PR^G_n$ l'ensemble des points fixes
r\'eguliers r\'epulsifs
de $G\circ F^n$.
Alors la suite de mesures 
$$\nu^+_n:=p_2^{-1}d_2^{-n}\sum_{z\in \PR^G_n}\delta_z$$
tend faiblement vers la mesure d'\'equilibre $\mu$ de $F$.
\end{theoreme}
\par
En prenant $G(z)=z$, on obtient le corollaire suivant qui
g\'en\'eralise un th\'eor\`eme de Lyubich \cite{Lyubich},
de Freire-Lopes-Ma\~n\'e
\cite{FreireLopesMane}
et de Briend-Duval \cite{BriendDuval1}. Ces
auteurs ont consid\'er\'e le cas des
endomorphismes holomorphes de $\P^k$.
\begin{corollaire} Soit $F$ une correspondance 
comme au th\'eor\`eme 2.1.
Notons $\PR_n$ l'ensemble des points p\'eriodiques r\'eguliers
r\'epulsifs de p\'eriode $n$ de $F$.
Alors la suite de mesures 
$$\nu_n:=d_2^{-n}\sum_{z\in \PR_n}\delta_z$$
tend faiblement vers la mesure d'\'equilibre $\mu$ de $F$.
\end{corollaire}
\par
Nous allons adapter
une m\'ethode de Briend-Duval \cite{BriendDuval1} qui traite
le cas des endomorphismes holomorphes de $\P^k$.
Il s'agit d'une application de la proposition
3.2 et du m\'elange sur une extension naturelle du
syst\`eme dynamique 
associ\'e \`a $F$ (voir \cite{CornfeldFominSinai}). Cependant, 
les propri\'et\'es d'invariance et 
de m\'elange de la mesure d'\'equilibre dans notre cas
sont plus faibles
que dans le cas des endomorphismes. Plus pr\'ecis\'ement, la
mesure $\mu$ n'est pas $F_*$-invariante, ce qui nous oblige \`a 
construire diff\'eremment la mesure associ\'ee \`a
l'extension naturelle de $F$. Cette mesure ne sera pas invariante
en g\'en\'eral.
\par
On a d\'efini dans le lemme 3.4 un polyn\^ome $P$ tel que tout
point de $(P\not=0)$ admet $d_2$ branches inverses
r\'eguli\`eres d'ordre 1.
Soit $Q$ un polyn\^ome analogue associ\'e \`a $G$. Notons 
$\Lambda$ l'op\'erateur de Perron-Frobenius associ\'e \`a $F$
et $\Lambda_G$ celui associ\'e \`a $G$. Notons $Z=\sum
Z_i^*$ la
$k$-cha\^{\i}ne holomorphe associ\'ee \`a $G$.
Notons $X_0$ le compl\'ementaire de
$G(\E')\cup (Q=0)$ o\`u $\E'$ est l'orbite de $(P=0)$ par $F$.
Du fait que $\mu$ ne charge pas les ensembles pluripolaires, 
$X_0$ est de $\mu$ mesure totale.
\par
Dans la suite, nous allons changer un peu
la terminologie "branche inverse r\'eguli\`ere" pour l'adapter \`a la suite 
$(G\circ F^n)$.
Soit $K_1$ un sous-ensemble connexe de $\C^k$.
On appelle {\it branche
inverse r\'eguli\`ere
d'ordre $n$} de $K_1$ toute suite $\B$ 
$$K_{-n},(\widehat{K}_{-n},i_n), K_{-n+1},
(\widehat{K}_{-n+1},i_{n-1}),
\ldots,K_{-1}, (\widehat{K}_{-1},i_1), K_0, (\widehat K_0,i), K_1$$
v\'erifiant les conditions suivantes
\begin{enumerate}
\item[(i)] Les ensembles 
$K_{-m}\subset\C^k$ sont connexes pour $m\geq -1$,
$\widehat{K}_{-m}\subset Y^*_{i_m}$ pour $m\geq 1$, $\widehat
K_0\subset Z_i^*$;
\item[(ii)]
$\pi_1$ d\'efinit une bijection de
$\widehat{K}_{-m}$ dans $K_{-m}$ et
$\pi_2$ d\'efinit une bijection de
$\widehat{K}_{-m}$ dans $K_{-m+1}$ 
pour $0\leq m\leq n$.
\end{enumerate}
Si $K_1$ n'est pas un ouvert (par exemple si $K_1$ est un point),
on exige aussi que les applications
ci-dessus d\'efinissent des bijections holomorphes
entre un voisinage de $\widehat K_{-m}$ et ses images.
Notons $\FF^{-m}_\B$ l'application
$$(\pi_{1|\widehat{K}_{-m}})\circ (\pi_{2|\widehat{K}_{-m}})^{-1}
\circ \cdots \circ
(\pi_{1|\widehat{K}_0})\circ (\pi_{2|\widehat{K}_0})^{-1}$$
pour $1\leq m\leq n$.
C'est une application holomorphe bijective de $K_1$ dans $K_{-m}$.
Les autres notations et la terminologie sont modifi\'ees de la
m\^eme mani\`ere.
\par
Si $K_1$ est une boule et si $K_{-n}$ est strictement contenu
dans $K_1$, alors l'application $\FF^{-n}_\B$ est contract\'ee
pour la m\'etrique de Kobayashi de $K_1$.
Elle admet donc un point fixe
attractif unique dans $K_{-n}$. C'est un point fixe r\'epulsif de 
la correspondance $G\circ F^n$.
Dans la suite, nous allons construire et compter les
branches inverses r\'eguli\`eres v\'erifiant cette propri\'et\'e.
Observons que tout point de $X_0$ admet exactement $p_2d_2^n$
branches inverses r\'eguli\`eres d'ordre $n$.
\par
Notons $X$ l'ensemble des branches
inverses r\'eguli\`eres
d'ordre infini d'un point $x_1\in X_0$
d\'efinies par les suites infinies 
$$\ldots, x_{-n}, (\widehat x_{-n},i_n),  x_{-n+1},
(\widehat x_{-n+1},i_{n-1}), \ldots,
 x_{-1}, (\widehat x_{-1},i_1), x_0, (\widehat x_0,i), x_1,$$
$x$ d\'esignera la suite infinie ci-dessus.
La notation $X_n$ d\'esignera
l'ensemble des branches inverses r\'eguli\`eres
$x^{(n)}$ d'ordre $n$ donn\'ees
par les suites finies
$$x_{-n}, (\widehat x_{-n},i_n),  x_{-n+1},
(\widehat x_{-n+1},i_{n-1}), \ldots,
 x_{-1}, (\widehat x_{-1},i_1), x_0, (\widehat x_0,i), x_1.$$
\par
Notons ${\cal A}_0$ la $\sigma$-alg\`ebre des bor\'eliens
de $X_0$. Consid\'erons la $\sigma$-alg\`ebre
${\cal A}$ de $X$ engendr\'ee
par les ensembles 
$$A_{-m}(S):=\{x\in X, x_{-m}\in S\}$$
o\`u $S$ est un \'el\'ement de ${\cal A}_0$ et $m\geq -1$.
D\'efinissons une mesure de probabilit\'e $\tilde\mu$ sur $X$.
Soient $X_{-m}$ des \'el\'ements de $\A_0$.  
Posons
$$A(S_1,S_0,\ldots,S_{-m}):=A_1(S_1)\cap\ldots\cap
A_{-m}(S_{-m})$$
et
$$\tilde\mu(A(S_1,S_0,\ldots, S_{-m})):=
\big\langle \mu,\ind_{S_1}\Lambda_G(\ind_{S_0}\Lambda(\ind_{S_{-1}}
(\Lambda \ind_{S_{-2}}\ldots 
(\Lambda \ind_{S_{-m}})\ldots)))\big\rangle.$$
Observons que la valeur de la fonction
$$p_2d_2^m\ind_{S_1}\Lambda_G(\ind_{S_0}\Lambda(\ind_{S_{-1}}
(\Lambda \ind_{S_{-2}}\ldots 
(\Lambda \ind_{S_{-m}})\ldots)))$$
en $x_1$ est \'egale au nombre
de branches inverses d'ordre $m$
$$x_{-m}, (\widehat x_{-m},i_m),  x_{-m+1},
(\widehat x_{-m+1},i_{m-1}), \ldots,
 x_{-1}, (\widehat x_{-1},i_1), x_0, (\widehat x_0,i), x_1$$
de $x_1$ qui v\'erifient $x_1\in S_1$, $x_0\in S_0$, $\ldots$,
$x_{-m}\in S_{-m}$. 
\par
On a la condition de compatibilit\'e suivante
\begin{eqnarray*}
\lefteqn{\tilde\mu(A(S_1,S_0\ldots,S_{-i+1},S_{-i}\sqcup S_{-i}',
S_{-i-1},\ldots, S_{-m})) =}\\
& = &\tilde\mu(A(S_1,S_0,\ldots,S_{-i+1},S_{-i},S_{-i-1},
\ldots, S_{-m})) + \\
& & +\tilde\mu(A(S_1,S_0,\ldots,
S_{-i+1},S_{-i}',S_{-i-1},\ldots, S_{-m}))
\end{eqnarray*}
lorsque $S_{-i}$, $S_{-i}'$ sont disjoints et $i\geq -1$.
D'apr\`es le th\'eor\`eme de consistence de Kolmogoroff,
$\tilde \mu$ s'\'etend, de mani\`ere unique,
en une mesure de probabilit\'e sur $X$.
\par
Notons $\Pi_n:X\longrightarrow X_n$ la projection
$\Pi_n(x):=x^{(n)}$ et $\tau_n:X_n\longrightarrow X_0$ la
projection $\tau_n(x^{(n)}):=x_1$.
L'application $\tau_n$ est de degr\'e $p_2d_2^n$ car on a
supprim\'e tous les points donnant naissance \`a de mauvaises branches
inverses.
Pour $\varphi$ une fonction sur $X_n$ posons
$$(\tau_n)_*\varphi(a):=\sum_{\tau_n(b)=a}\varphi(b).$$
\begin{lemme}
Soit $A$ un \'el\'ement de $\A$. Alors on a
$$\tilde\mu(A)=\lim_{n\rightarrow\infty} \int p_2^{-1}
d_2^{-n}(\tau_n)_*\ind_{\Pi_n(A)}\d\mu.$$
\end{lemme}
\begin{preuve} Il suffit de consid\'erer le cas o\`u $A$ est du type 
$A(S_1,S_0,\ldots,S_{-m})$. On a alors 
$$\tilde\mu(A)=\int p_2^{-1}
d_2^{-n}(\tau_n)_*\ind_{\Pi_n(A)}\d\mu$$
pour tout $n\geq m$. D'o\`u le lemme.
\end{preuve}
\begin{lemme} Soit $S$ un \'el\'ement de $\A_0$ et
soit $A$ un \'el\'ement de $\A$. Alors on a
$\lim_{n\rightarrow\infty} \tilde\mu(A_{-n}(S)\cap
A)=\mu(S)\tilde\mu(A)$.
\end{lemme}
\begin{preuve} Il suffit de consid\'erer le cas o\`u
$A=A(S_1,S_0,\ldots, S_{-m})$.
Posons $c:=\mu(S)$ et $\varphi:=\ind_S-c$.
On a pour tout $n\geq m$
\begin{eqnarray*}
\lefteqn{|\tilde\mu(A_{-n}(S)\cap A) - c\tilde\mu(A)|=}\\
& = & \Big|\int\ind_{S_1}\Lambda_G(\ind_{S_0}\Lambda(\ind_{S_{-1}}
\ldots\Lambda(\ind_{S_{-m}}\Lambda^{n-m}\ind_S)\ldots))\d\mu-\\
& &
-c\int\ind_{S_1}\Lambda_G(\ind_{S_0}\Lambda(\ind_{S_{-1}}\ldots
\Lambda(\ind_{S_{-m}})\ldots))\d\mu \Big|\\
& = & \Big|\int\ind_{S_1}\Lambda_G(\ind_{S_0}\Lambda(\ind_{S_{-1}}
\ldots\Lambda(\ind_{S_{-m}}\Lambda^{n-m}\varphi)\ldots))
\d\mu\Big|\\
&\leq & \int \Lambda_G\Lambda^m|\Lambda^{n-m}\varphi| \d\mu.
\end{eqnarray*}
Notons $\C^\perp$ l'orthogonal de $\C$ dans $\Ltwo(\mu)$.
On a $\varphi\in\C^\perp$, $\|\varphi\|_\infty\leq 1+|c|$ et donc
$0\leq
\Lambda_G\Lambda^m|\Lambda^{n-m}\varphi|\leq 1+|c|$. 
D'apr\`es le th\'eor\`eme de convergence domin\'ee,
il suffit de montrer que pour tout $\psi\in\C^\perp$, $\Lambda^n\psi$ tend vers $0$
$\mu$-presque partout (extraire des sous-suites si n\'ecessaire). 
Or, d'apr\`es le th\'eor\`eme 2.10,
ceci est
vrai pour $\psi$ lisse. Du fait que le sous-espace des fonctions
lisses est dense dans $\C^\perp$, il suffit de montrer que 
$\|\Lambda\|_{\Ltwo(\mu)}\leq 1$. Gr\^ace \`a l'in\'egalit\'e de
Cauchy-Schwarz et \`a l'invariance de $\mu$, on a
\begin{eqnarray*}
\int |\Lambda\psi(z)|^2\d\mu(z) & = & \int \Big|\frac{1}{d_2}
\sum_{w\in F^{-1}(z)}\psi(w)\Big|^2\d\mu(z) \\
& \leq & \int \frac{1}{d_2}\Big(\sum_{w\in
F^{-1}(z)}|\psi(w)|^2 \Big)\d\mu(z)\\
& = & \int\Lambda (|\psi|^2) \d\mu = \int |\psi|^2\d\mu.
\end{eqnarray*}
Ceci termine la preuve du lemme.
\end{preuve}
{\bf Fin de la d\'emonstration du th\'eor\`eme 4.2.}
Posons, pour tout $\delta>0$, $E_\delta$
l'ensemble des $x\in X$ donn\'es
par les suites
$$\ldots, x_{-n}, (\widehat x_{-n},i_n),  x_{-n+1},
(\widehat x_{-n+1},i_{n-1}), \ldots, x_{-1},
(\widehat x_{-1},i_1), x_0, (\widehat
x_0,i), x_1$$
tels que la boule $\Ball(x_1,\delta)$ admette une branche inverse
r\'eguli\`ere
d'ordre infini, de taille $(l^{-m/4})$ et accroch\'ee \`a $x$.
D'apr\`es la proposition 3.6,
l'union $\cup_{\delta>0} E_\delta$ est de
$\tilde\mu$ mesure totale. Soit $\pi:X\longrightarrow X_0$ la projection
$\pi(x):=x_1$. Posons
$\mu_\delta:=\pi_*(\ind_{E_\delta}\tilde\mu)$. Cette
mesure tend vers $\mu$ quand $\delta$ tend vers $0$. D'apr\`es la
proposition 4.1, quand $n$ est assez grand,
les mesures $\nu^+_n$ sont de masse au plus 1.
Pour prouver le th\'eor\`eme 4.2, on montre que toute
valeur adh\'erente $\nu^+$ \`a la suite $(\nu^+_n)$ v\'erifie 
$\mu_\delta\leq \nu^+$. Fixons $\delta$, $0<\delta<1$,
$\epsilon$, $0<\epsilon<\delta$  et un point $z_1\in X_0$.
Soit $r$ v\'erifiant $0<2r<\delta-\epsilon$.
Il suffit de montrer que
$(1-\epsilon)\mu_\delta(\overline \Ball(z_1,r)) \leq
\nu^+_n(\Ball(z_1,r+\epsilon))$ pour tout $n$ assez grand.
\par
Posons $S:=\overline \Ball(z_1,r)$ et $A:=\pi^{-1}(S)\cap
E_\delta$. D'apr\`es le lemme 4.5, on a
$$\lim_{n\rightarrow \infty} \tilde\mu (A_{-n}(S)\cap A)=
\mu(S)\tilde\mu(A)=\mu(S)\mu_\delta(S).$$
Soit $n$ assez grand tel que $l^{-n/4}<\epsilon$ et
$\tilde\mu(A_{-n}(S)\cap A)>
(1-\epsilon) \mu(S)\mu_\delta(S)$. Notons pour tout $z\in X_0$,
$s(z)$ le nombre de branches inverses r\'eguli\`eres $\B$ d'ordre $n$,
de taille $(l^{-{m/4}})$ de $\Ball(z,\delta)$ telles que
$\FF^{-n}_\B(z)\in S$. On d\'eduit de la derni\`ere in\'egalit\'e
et du lemme 4.4 que pour $n$ assez grand
$$\int p_2^{-1}d_2^{-n}s(z)\ind_S \d\mu \geq
(1-\epsilon)\mu(S)\mu_\delta(S).$$
Par cons\'equent, il existe un point $x_1\in S$ tel que
$s(x_1)\geq (1-\epsilon) p_2d_2^n\mu_\delta(S)$. Il existe donc
au moins $(1-\epsilon)p_2d_2^n\mu_\delta(S)$ branches
inverses r\'eguli\`eres $\B$ d'ordre $n$,
de taille $(l^{-{m/4}})$ de $\Ball(x_1,\delta)$ telles que
$z_1:=\FF^{-n}_\B(x_1)\in S$. L'application $\FF^{-n}_\B$ d\'efinit une
application bijective de $\Ball(x_1,\delta)$ dans
$W:=\FF^{-n}_\B \Ball(x_1,\delta)$. Or $W$ est un ensemble
de diam\`etre plus petit que $\epsilon$ et contient le point
$z_1\in S$. On en
d\'eduit que $\overline W$ est contenu dans
$\Ball(z_1,r+\epsilon) \subset \Ball(x_1,\delta)$.
Par cons\'equent, $\FF^{-n}_\B$ admet un point fixe attractif
unique dans $\Ball(z_1,r+\epsilon)$.
On en d\'eduit que
$\nu^+_n(\Ball(z_1,r+\epsilon))\geq (1-\epsilon)
\mu_\delta(\overline\Ball(z_1,r))$. La d\'emonstration du
th\'eor\`eme 4.2 est achev\'ee.
\par
\hfill $\square$
\par
\begin{remarque} \rm
Dans le corollaire 4.3, on peut remplacer l'ensemble $\PR_n$ par
l'ensemble $\PR_n\cap\supp(\mu)$. En effet, dans la preuve du
th\'eor\`eme 4.2 pour le cas o\`u $G(z)=z$, il suffit de
consid\'erer le point $x_1$ appartenant \`a $\supp(\mu)$. Or
$\supp(\mu)$ est invariant par $F^{-1}$. Donc
les points fixes obtenus
appartiennent n\'ecessairement \`a $\supp(\mu)$.
On peut aussi remplacer $\PR_n$ par l'ensemble $\PR'_n$ des
point $x\in\PR_n$ dont la p\'eriode minimale est \'egale \`a
$n$. En effet, d'apr\`es la proposition 4.1, on a
$\#(\PR_n\setminus \PR_n') =\o(d_2^n)$. 
\end{remarque}
\par
Les corollaires suivants donnent des interpr\'etations
g\'eom\'etriques des r\'esultats obtenus.
\begin{corollaire} Soient $F$ une correspondance comme au
th\'eor\`eme 2.1 et $G$ une autre correspondance polynomiale de
degr\'e topologique $(p_1,p_2)$ associ\'ee \`a une
$k$-cha\^{\i}ne holomorphe $Z$. Notons $Y_n$
la $k$-cha\^{\i}ne holomorphe associ\'ee \`a $F^n$ et $R_n$
l'intersection de $Y_n$ avec $Z$. Alors la suite de mesures
$$p_1^{-1}d_2^{-n}\sum_{z\in R_n}\delta_{\pi_1(z)}$$
tend faiblement vers la mesure d'\'equilibre $\mu$ de $F$.
\end{corollaire}
\begin{preuve} Soit $\overline G$ la correspondance adjointe de
$G$. 
Chaque point fixe de $\overline G\circ F^n$ est associ\'e
\`a un couple de points
$(x,y)\in Y_n$ et
$(y,x)\in \overline Z$. On peut
donc l'associer au point $(x,y)$ dans l'intersection de $Y_n$
avec $Z$. Il suffit maintenant
d'appliquer le th\'eor\`eme 4.2 en rempla\c
cant $G$ par $\overline G$.  
\end{preuve}
\begin{corollaire} Soit $F$ une correspondance comme au
th\'eor\`eme 2.1.
Soient $Y_n$ la $k$-cha\^{\i}ne holomorphe associ\'ee \`a $F^n$
et $[Y_n]$ le courant d'int\'egration sur
$Y_n$. 
Alors le courant $d_2^{-n}[Y_n]$ tend faiblement vers le courant
$\pi_1^*(\mu)$.
\end{corollaire}
\begin{preuve} Pour toute matrice $A$ inversible de norme petite et
de rang $k$,
notons $\pi_A:\C^k\times\C^k\longrightarrow \C^k$ l'application d\'efinie par
$\pi_A(z^1,z^2):=z^2-Az^1$.
Soit $\alpha:=\d z_1\wedge \d\overline z_k
\wedge \ldots \wedge \d z_k \wedge \d \overline z_k$. Posons
$\alpha_A:=(\pi_A)^*\alpha$. Les formes $\Phi\alpha_A$ avec
$\Phi$ une fonction \`a support compact, engendrent l'espace des
$(k,k)$-formes \`a support compact dans
$\C^k\times\C^k$. Il suffit de montrer que
$$\lim_{n\rightarrow\infty} \langle d_2^{-n} [Y_n]-\pi_1^*(\mu),
\Phi\alpha_A\rangle =  0.$$
Notons $E_z:=\pi_A^{-1}(z)$ pour tout $z\in\C^k$ et $\Phi_z$ la
restriction de $\Phi$ \`a $E_z$.
Notons \'egalement $[Y_n^z]$ et $\nu^z$ les mesures obtenues
comme intersections de $[Y_n]$ et de $\pi_1^*(\mu)$ avec $E_z$. 
D'apr\`es le th\'eor\`eme de
Fubini on a
$$\langle d_2^{-n} [Y_n]-\pi_1^*(\mu), \Phi\alpha_A\rangle =
\int \langle d_2^{-n}[Y_n^z]-\nu^z,\Phi_z \rangle \d z_1 \wedge \d
\overline z_1 \wedge \ldots \wedge \d z_k
\wedge \d \overline z_k.$$
D'apr\`es le corollaire 4.7, $\langle d_2^{-n}[Y_n^z]-\nu^z,\Phi_z
\rangle$ tend vers $0$ pour tout $z$. Le th\'eor\`eme de
convergence domin\'ee implique que 
$$\lim_{n\rightarrow\infty}
\langle d_2^{-n} [Y_n]-\pi_1^*(\mu), \Phi\alpha_A\rangle
=0.$$ 
\end{preuve}
\section{Ensembles d'unicit\'e pour les polyn\^omes}
Soit $P$ un polyn\^ome de degr\'e $\deg(P)\geq 2$.
Soient $K_P$ l'ensemble des points
d'orbite born\'ee de $P$ et $J_P$ son bord topologique. L'ensemble
$K_P$ s'appelle {\it ensemble de Julia rempli} de $P$ et $J_P$
est {\it l'ensemble de Julia } de $P$. D'apr\`es la th\'eorie de
Fatou, $J_P$ est le plus petit compact qui
contient plus d'un point et qui v\'erifie $P^{-1}(J_P)=J_P$.
L'ensemble $K_P$ est le plus grand compact v\'erifiant
$P^{-1}(K_P)=K_P$.
On dira qu'un compact $K$ est {\it un ensemble d'allure Julia} de $P$
si $K$ contient au moins deux points et si $P^{-1}(K)=K$.
Il est clair qu'un tel ensemble $K$ n'est pas
un ensemble d'unicit\'e et v\'erifie 
$J_P\subset K\subset K_P$.
\par
Consid\'erons deux polyn\^omes distincts non constants
$f$, $g$ avec $\deg
(f)\geq \deg (g)$ et un compact $K$, $\#K\geq 2$, 
v\'erifiant $f^{-1}(K)=g^{-1}(K)$.
Dans \cite{Ostrovskii}
Ostrovskii, Pakovitch et Zaidenberg ont montr\'e que si
$\deg(f)=\deg(g)$ il existe une rotation $R$ pr\'eservant $K$
telle que $f=R\circ g$. 
Si la capacit\'e
logarithmique de $K$ est strictement positive,
nous avons
montr\'e \cite{Dinh}
qu'il existe un polyn\^ome $P$ tel que $f=P\circ g$ et
$P^{-1}(K)=K$ sauf pour les deux cas exceptionnels ci-dessous. D\'esignons par
$z$, $Q$, $d$, $d'$ et $a$
une coordonn\'ee, un polyn\^ome, deux entiers naturels
et un nombre complexe convenables
\par
{\bf Cas 1}.
$K$ est une r\'eunion de cercles de centre $0$ et $f(z)=Q(z)^d$,
 $g=aQ(z)^{d'}$.
\par
{\bf Cas 2}. $K$ est le segment $[-1,1]$ et $f=\pm\T_d\circ Q$,
$g=\pm\T_{d'}\circ Q$
o\`u $\T_m$ est le polyn\^ome de Tchebychev
d\'efini par la relation $\T_m(\cos z)=\cos(mz)$.
\par
On voit que hors de ces deux cas
si $\deg(f)>\deg(g)$, $K$
est un ensemble d'allure Julia du polyn\^ome $P$.
\par
On d\'eduit de ces r\'esultats le th\'eor\`eme suivant.
\begin{theoreme} Soient $f$, $g$ deux polyn\^omes  
non constants avec $\deg(f)\geq \deg(g)$.
Soit $K$ un ensemble compact infini de $\C$ v\'erifiant
$f^{-1}(K)=g^{-1}(K)$. Alors il existe un polyn\^ome $P$ tel que
$f=P\circ g$ et $P^{-1}(K)=K$ sauf si $(f,g,K)$ appartient \`a l'une
des deux classes d\'ecrites ci-dessus.
\end{theoreme}
\begin{preuve} Si $\deg(f)=\deg(g)$, ce th\'eor\`eme se r\'eduit
au r\'esultat d'Ostrovskii-Pakovitch-Zaidenberg.
Supposons que $\deg(f)>\deg(g)$. Posons $F:=f\circ g^{-1}$. C'est
une correspondance polynomiale
dont l'exposant de Lojasiewicz est \'egal
\`a $\deg(f)\deg(g)^{-1}$. Elle v\'erifie l'hypoth\`ese du
corollaire 3.10. On a $F^{-1}(K)=K$ et $F(K)=K$. Comme $K$ est
un ensemble infini, il n'est pas contenu dans
l'ensemble $\E_0$.
D'apr\`es le corollaire 3.10 et la remarque 3.12,
les relations
$F^{-1}(K)=F(K)=K$ impliquent que $K$ contient le support de la mesure
d'\'equilibre $\mu$ de $F$.
La mesure $\mu$ \'etant PB,
la capacit\'e logarithmique de son support
est strictement positive. Les
r\'esultats cit\'es ci-dessus permettent de conclure.
\end{preuve}
\begin{corollaire} Un sous-ensemble compact infini de $\C$ est un
ensemble d'unicit\'e si et seulement s'il n'est pas un ensemble
d'allure Julia et s'il n'est invariant par aucune rotation. 
\end{corollaire}

Tien-Cuong Dinh, B\^at. 425 -- Math\'ematique, Universit\'e
Paris-Sud, 91405 Orsay, France. E-mail:
Tiencuong.Dinh@math.u-psud.fr 
\end{document}